\documentclass[11pt,reqno]{amsart}
\usepackage[dvips]{graphicx,color}
\usepackage{amssymb,amscd,latexsym,epsfig}
\usepackage[arrow,matrix,line]{xy}
\usepackage[mathscr]{euscript}

\DeclareFontFamily{OT1}{rsfs}{}
\DeclareFontShape{OT1}{rsfs}{n}{it}{<-> rsfs10}{}
\DeclareMathAlphabet{\curly}{OT1}{rsfs}{n}{it}

\newcommand{\rt}[1]{\stackrel{#1\,}{\rightarrow}}
\newcommand{\Rt}[1]{\stackrel{#1\,}{\longrightarrow}}
\newcommand\comp{{\!\,}_{{}^\circ}}
\newcommand\B{\operatorname{Bl}}
\newcommand\C{\mathbb C}
\newcommand\Q{\mathbb Q}
\newcommand\R{\mathbb R}
\newcommand\G{\mathbb G}
\renewcommand\S{\mathscr S}
\newcommand\Z{\mathcal Z}
\newcommand\PP{\mathbb P}
\renewcommand\P{\mathcal P}

\newcommand\X{\curly X}
\newcommand\XX{\,\widehat{\!\curly X}{}}
\newcommand\Y{\curly Y}
\renewcommand\L{\mathcal L}

\renewcommand\O{\mathcal O}
\newcommand\I{\curly I}
\newcommand\into{\hookrightarrow}
\newcommand\res{\arrowvert_}

\renewcommand\_{^{\ }_}

\newcommand{\rk}{\operatorname{rank}}

\newcommand{\im}{\operatorname{im}}

\newcommand{\Proj}{\operatorname{Proj}\,}
\newcommand{\Spec}{\operatorname{Spec}\,}
\newcommand{\Hilb}{\operatorname{Hilb}}
\newcommand{\Grass}{\operatorname{Grass}}

\makeatletter \@addtoreset{equation}{section} \makeatother

\newtheorem{thm}[equation]{Theorem}
\newtheorem{lem}[equation]{Lemma}
\newtheorem{cor}[equation]{Corollary}
\newtheorem{prop}[equation]{Proposition}

\theoremstyle{definition}
\newtheorem{defn}[equation]{Definition}
\newtheorem{rmk}[equation]{Remark}
\newtheorem{rmks}[equation]{Remarks}

\title[The Hilbert-Mumford criterion]
{\textbf{A study of the Hilbert-Mumford criterion for the stability of
projective varieties}}
\author[J. Ross and R. P. Thomas]{Julius Ross and Richard Thomas}

\begin{document}

\begin{abstract} \noindent
We make a systematic study of the Hilbert-Mumford criterion for different
notions of stability for polarised algebraic
varieties $(X,L)$; in particular for K- and Chow stability. For each type of stability this leads to a concept of
slope $\mu$ for varieties and
their subschemes; if $(X,L)$ is semistable then $\mu(Z)\le\mu(X)$ for all
$Z\subset X$. We give examples such as curves, canonical models and Calabi-Yaus.
We prove various foundational technical results towards understanding
the converse, leading to partial results; in particular this gives a geometric
(rather than combinatorial) proof of the stability of smooth curves.
\end{abstract}

\maketitle


\section{Introduction; slope stability}

Geometric Invariant Theory \cite{GIT} has been very successful in forming
moduli spaces of (semi)stable coherent sheaves over polarised
algebraic varieties $(X,L)$ \cite{HL}. Moduli space is constructed as a quotient
of a subset of a Quot scheme by a group, and the Hilbert-Mumford criterion
is applied to 1-parameter subgroups. Their weights are found to be dominated
by positive linear combinations of weights of particularly simple 1-parameter
subgroups corresponding to a degeneration
of a sheaf $E$ into a splitting $F\oplus E/F$, for some subsheaf $F\le E$.
Thus stability of sheaves is governed by subsheaves; calculating the corresponding
weights leads to the notion of slope stability (the exact form of the slope
depending on the linearisation used on Quot).

Forming moduli spaces of (semi)stable varieties themselves using GIT has
proved much more difficult, and has mainly been accomplished for canonically
polarised varieties using the Chow linearisation, due to work of Mumford
\cite{GIT, Mu}, Gieseker \cite{Gi}, and, in a little more generality, Viehweg
\cite{V}. Roughly
speaking one expects varieties polarised by their canonical bundle $L=K_X$
to be automatically stable since their moduli functor is
already separated because of the birational invariance of spaces of sections
of powers of the canonical bundle. In the general case no geometric
criterion for (in)stability has emerged. This is because the Hilbert-Mumford
criterion has not been successfully simplified or interpreted for varieties;
instead Viehweg proved deep positivity results to produce the group-invariant
sections of the appropriate line bundle directly (for varieties with semi-ample
canonical bundle). Kollar \cite{Ko} and others have turned to other methods
for producing moduli of varieties, but new impetus to understanding stability
has come from the link between K-stability and the existence
of K\"ahler-Einstein and constant scalar curvature K\"ahler metrics \cite{Ti1,
Do1} to which we apply our methods in \cite{RT}.

Our approach uses the Hilbert-Mumford criterion, as pioneered by
Mumford \cite{Mu}. He calculates the relevant weights in terms of 
blow-ups of $X\times\C$ in subschemes supported on thickenings of $X\times\{0\}$;
one of our main results (Corollary \ref{ref}) reduces the calculation to
a sum of weights of blow-ups \emph{in the scheme-theoretic central fibre}.
(In a sense which is made clear by the proof the construction
turns the horizontal thickenings of Mumford ``vertical", into the
central fibre.)

Taking just one such blow-up in a subscheme $Z\subset X$ gives the ``deformation
to the normal cone of $Z$", analogous
to the simple 1-parameter subgroups that arise in the GIT of sheaves.
This gives a numerical condition for $Z$ to destabilise $X$, and so
a notion of ``slope stability", by analogy with the sheaf theory which we
now review briefly.

For a sheaf $E$ over $(X,L)$, the \emph{reduced} Hilbert polynomial $p\_E$
is the monic version of the Hilbert polynomial
$\P\_{\!E}(r)=\chi(E\otimes L^r)=a_0r^n+a_1r^{n-1}+\ldots\,$:
$$
p\_E(r)=\frac{\chi(E\otimes L^r)}{a_0}=r^n+\mu(E)r^{n-1}+\ldots\,,
$$
and the \emph{slope} of $E$ is its leading nontrivial coefficient,
\begin{equation} \label{mu}
\mu(E)=a_1/a_0.
\end{equation}
(This differs from the usual definition deg$\,(E)/\rk(E)$ by unimportant
terms that depend only on the geometry of $(X,L)$.)
Then we say that $E$ is \emph{semistable} if for all proper coherent subsheaves
$F\le E$,
$$
p\_F\preceq p\_E\,.
$$
For \emph{Gieseker semistability}, $\preceq$ means there exists an $r_0>0$
such that $p\_F(r)\le p\_{E}(r)$ for all $r\ge r_0$; for \emph{slope semistability}
the inequality is at the level of the $r^{n-1}$ coefficients:
\begin{equation} \label{mu2}
\mu(F)\le\mu(E).
\end{equation}
If the inequalities are all strict then $E$ is \emph{stable}, and this agrees
with the appropriate GIT notions for different choices of linearisation (i.e.\
choice of equivariant line bundle on the Quot scheme; in fact Jun Li's
line bundle \cite{Li} is only semi-ample, but he extends GIT to this case.)

Due to different choices
of linearisation and asymptotics there are also many notions of stability
(Hilbert, Chow, asymptotic Hilbert, asymptotic Chow and K-stability) for
varieties, defined in the next section, and so many versions of slope. We
first describe the one relevant to K-stability.

Given a subscheme $Z$ of a polarised algebraic variety $(X,L)$, define
$$
\epsilon(Z)=\sup\{c\in\Q_{\,>0}\colon L^r\otimes\I_{\!Z}^{cr}\text{ is globally
generated } \forall r\gg0 \text{ with } cr\in\mathbb N\}.
$$
The Hilbert-Samuel polynomial for \emph{fixed} $x$,
\begin{equation} \label{HS}
\chi(L^r\otimes\I^{xr}_{\!Z})=a_0(x)r^n+a_1(x)r^{n-1}+\ldots\,,\quad r\gg0,
\ rx\in\mathbb N,
\end{equation}
defines $a_i(x)$ which are polynomials in $x$ (see Section \ref{K})
and so extend by the same formulae to
$x\in\R$. Then analogously to (\ref{mu}) we define the \emph{K-slope} of
$\I_Z$ (with respect to $L$ and $c\in(0,\epsilon(Z)]$) to be
\begin{equation} \label{def:slope}
\mu_c(\I_Z)=\mu_c(\I\_Z,L):=\frac{\int_0^c\big(a_1(x)+\frac{a_0'(x)}2\big)dx}
{\int_0^ca_0(x)dx}\,.
\end{equation}
Setting $Z=\emptyset$ defines the slope of $X$ (precisely: of $\O_X$ with
respect to $L$ and $c$) as
\begin{equation} \label{def:X}
\mu(X)=\mu(X,L)=\frac{a_1}{a_0}\,,
\end{equation}
independently of $c$. Here $a_i$ are the coefficients of the
Hilbert polynomial $\chi(L^r)=a_0r^n+a_1r^{n-1}+\ldots\,,\ a_0=a_0(0)$
and, for $X$ normal, $a_1=a_1(0)$ (\ref{normal}).

Setting $\mu(Z):=\sup\_{0<c\le\epsilon(Z)}\mu_c(\I_Z)$, we say that
$(X,L)$ is \emph{K-slope semistable} if for all proper $Z\subset X$,
$$
\mu(Z)\le\mu(X), \qquad\qquad\text{i.e.}\quad \frac{\int_0^c\left(a_1(x)+
\frac{a_0'(x)}2\right)dx}{\int_0^ca_0(x)dx}\le\frac{a_1}{a_0}\quad
\forall c\in(0,\epsilon(Z)].
$$
Cf. (\ref{mu2}).
The definition of slope stability is slightly trickier (Definition \ref{defstab});
for this reason we work with $\mu_c(\I_Z)$ instead of $\mu(Z)$. Then (Theorem
\ref{thm:kstableslopestable}) $X$ is slope (semi)stable if it is K-(semi)stable.

The $a_0'/2$ ``correction term" in the definition of slope arises from the
difference between the Hilbert polynomial of a 2-component normal crossing
variety and the sum of those of its components (whereas for sheaves, $\P\_{\!E}=
\P\_{\!F}+\P\_{\!E/F}$ for any $F\le E$). Simon Donaldson
pointed out that his analysis of stability of toric varieties \cite{Do2}
throws up a boundary term similar to our $(a_0(c)-a_0)/2=\int_0^ca_0'/2$;
we explain this in \cite{RT}.

Similarly for $X\subseteq\PP^N=\PP(H^0(X,L)^*)$ (embedded in projective space
by the space of sections of its polarisation $L=\O_X(1)$ for ease of exposition;
see Section \ref{sec:chow} for the general case), a subscheme $Z\subset
X$ and an \emph{integer} $0<c\le\epsilon(Z)$,
we define the \emph{Chow slope} of $\I_Z$ to be
$$
Ch_c(\I_Z):=\frac{\sum_{i=1}^ch^0(\I_{\!Z}^i(1))}{\int_0^ca_0(x)dx}\,.
$$
Setting $Z=\emptyset$ gives $Ch(\O_X)=Ch(X):=\frac{h^0(\O_X(1))}{a_0}$;
then Chow (semi)stability implies Chow slope (semi)stability (Theorem
\ref{thmchow}):
$$
Ch_c(\I_Z){\ }^{\ <\,}_{(\le)}\ Ch(X).
$$
Asymptotic Chow is more like Gieseker stability and so our slope criterion
for that (\ref{chow}) is more complicated.

Section \ref{K} describes the various stability notions uniformly, while
Section \ref{Kslope} introduces slope stability via the deformation to the
normal cone. Arbitrary 1-parameter subgroups are studied in Section \ref{testconfigs},
calculating their weights in terms of a sequence of simple blow ups in Corollary
\ref{ref}. In Section \ref{converse} these weights are shown to be those
of a deformation to the normal cone under certain circumstances, giving a
partial
converse to ``stability $\Rightarrow$ slope stability". We used to think
this could be done in general, but the failure of the \emph{thickenings}
of certain flat families of subschemes to themselves be flat prevents  
us from carrying out the programme in full. We study
when the thickenings \emph{are} flat, and get round the problem with a basechange
trick in some situations. This deals with the curve case completely, i.e.\
stability and slope stability are equivalent there,
giving geometric proofs of the K- and asymptotic Chow stability of curves
of genus $\ge1$. As far as we know all previous proofs used analysis and
combinatorics
respectively. Section \ref{sec:chow} is devoted to Chow stability
and Section \ref{egs} to examples. Many more examples, such as projective
bundles, appear in \cite{RT}, in particular showing that K- and slope
stability are also equivalent for projective bundles over a curve.

\noindent \textbf{Acknowledgements.} 
We would like to thank Brian
Conrad, Kevin Costello, Dale Cutkosky, Simon Donaldson, David Eisenbud,
Daniel Huybrechts, Frances Kirwan, Miles Reid and
Bal\'azs Szendr\"oi for useful conversations. The book
\cite{HL} and paper \cite{Mu} have been very useful to us. The authors
were supported by an EPSRC PhD studentship and a Royal Society
university research fellowship respectively. The second author would also
like to thank Oscar Garcia-Prada and Peter Newstead for supporting a visit
to CSIC, Madrid, where much of this work was done.

\section{Notation} \label{notation}

Throughout this paper $Z$ will denote a closed subscheme of a proper
irreducible polarised scheme $(X,L)$ of dimension $n=\dim X$; for much of
the paper $X$ will be a normal irreducible variety. By $jZ$ we mean the subscheme
which is the $j$-fold thickening of $Z$ defined by its ideal sheaf $\I_{jZ}:=
\I_{\!Z}^j$. We denote
the blow up along $Z$ by $\pi\colon\widehat X\to X$, with exceptional divisor
$E$. Then $\pi^*\I_{\!Z}^j=\O(-jE)$ and there exists a $p$ such that
$\pi_*\O(-jE)=\I_{\!Z}^j$ for all $j\ge p$.  For convenience
we often suppress pullback maps, mix multiplicative and additive
notation for line bundles, and use the same letter to denote a divisor
and the associated line bundle.  For example on $\widehat{X}$, we
denote $(\pi^*L(-E))^{\otimes k}$ by $L^k-kE$. Worse still, where it
does not cause confusion, $L^k$ may also denote $c_1(L)^{\cap k}$.

For brevity we often denote sheaf cohomology on a space $X$ by $H^i_X$;
this never means local cohomology.

Any vector space $V$ with a $\C^\times$-action splits into one-dimensional
weight spaces $V=\bigoplus_iV_i$, where $t\in\C^\times$ acts on $V_i$ by
$t^{w_i}$. The $w_i$ are the \emph{weights} of the action, and $w(V)=\sum_iw_i$
is the \emph{total weight} of the action; i.e.\ the weight of the induced
action on $\Lambda^{\max}V$.

On any family over $\C=\Spec\C[t]$, $t$ will denote the pullback of the
standard coordinate on $\C$.

A line bundle $L$ is \emph{semi-ample} \cite{La} if its high
powers are globally generated (i.e.\ basepoint free). In this paper all
semi-ample line bundles will have the additional property that
the contraction they define is birational; that is $L$ is the pull back of
an ample line bundle from a birational scheme. (It is important that this
contraction can be trivial, i.e.\ semi-ample includes ample.)

Given a (semi-)ample line bundle $L$ on $X$, the \emph{Seshadri constant}
of $Z$ is
\begin{eqnarray} \nonumber
 \epsilon(Z) = \epsilon(Z,X,L)&=&\sup\,\{c\in\Q_{\,>0}\colon L^r\otimes
 \I_{\!Z}^{cr} \text{ is globally generated for } r\gg0\} \\ \label{def:seshadri}
&=&\max\,\{ c \in \Q_{\,>0}: L - cE \,\, \text{is nef on }\widehat X \}.
\end{eqnarray}

Given a pair of ideals $J\subset I\subset\O_X$, we say that $J$ \emph{saturates}
I if there exists $i>0$ such that $JI^{i-1}=I^i$. Equivalently, 
on the blow up $p\colon\widehat X\to X$ of $X$ along $I$ with exceptional
divisor $E$, the natural inclusion $p^*J\to\O(-E)$ is an isomorphism.
(This equivalence is a tautology: the zero set of the sections
of $p^*J\subseteq\O(-E)$ has homogeneous ideal the saturation of
$\oplus_iJI^{i-1}\subseteq\oplus_iI^i$ (\cite{Ha} Exercise II.5.10); this
is all of $\oplus_iI^i$ if and only if the zero set is empty, if and only
if the sections generate $\O(-E)$.)

Similarly, give a line bundle $L$ on $X$, the global sections of $L\otimes
I$ generate a subsheaf $L\otimes J\subset L\otimes I$, defining an ideal
$J$. We say that \emph{the global sections of $L\otimes I$
saturate} $I$ if $J$ saturates $I$. Equivalently, the global
sections of $L\otimes I$ generate the line bundle $L(-E)$ on $\widehat
X$. This is weaker than (i.e.\ is implied by) $L\otimes I$ being globally
generated.

\section{Hilbert, Chow and K-stability, and test configurations} \label{K}

Fix a polynomial $\P$ of degree $n$ and consider any $n$-dimensional
proper polarised scheme $(X,L)$ whose Hilbert polynomial equals $\P$:
$$
\P(r)=\chi\_X(L^r) = a_0r^n+a_1r^{n-1}+\ldots\,,
$$
where
$a_0=\frac{1}{n!}\int_X c_1(L)^n=\frac{L^n}{n!}$ and, for smooth $X$,
$$
a_1=\frac{1}{2(n-1)!}\int_X c_1(X)c_1(L)^{n-1}=
-\frac{K_X.L^{n-1}}{2(n-1)!}\,.
$$
Fix $r>0$ such that $L^r$ is both very ample on $X$, and 
\begin{equation} \label{rlarge}
H^i(L^r)=0 \quad\text{for }i>0, \quad\Rightarrow\ H^0(L^r)^*\cong\C^{\P(r)}.
\end{equation}
Then $L^r$ embeds $X$ in $\PP=\PP^{\P(r)-1}$, defining a point of the
Hilbert scheme $\Hilb=\Hilb_{\P',K}$ of subschemes of $\PP$ with
Hilbert polynomial $\P'(K)=\P(Kr)$.  Then there is a $K_0$ such that for
\emph{all} points $\{X\}$ of Hilb and $K\ge K_0$ we have an exact sequence
\begin{equation} \label{eq:seshilb}
0\to H^0_\PP(\I_X(K))\to S^K \C^{\P(r)*}\cong S^K H_X^0(L^r)\to
H^0_X(L^{Kr})\to0.
\end{equation}
This (see for example \cite{V}) defines Hilb as a closed subscheme of the
Grassmannian 
$$
\G=\Grass(S^K \C^{\P(r)*},\P(Kr)).
$$

So $(X,L^r)$ and a choice of isomorphism $H^0(L^r)^*\cong\C^{\P(r)}$
give us a point $x=x_{r,K}$ in $\G$, with different choices of isomorphism
corresponding (up to scale) to the orbits of
$SL(\P(r),\C)$ on $\PP^{\P(r)-1}$. A $g\in SL(\P(r),\C)$ acting on $\C^{\P(r)}$
induces an action
\begin{equation} \label{convention}
(S^Kg^*)^{-1} \quad\text{on }\ S^K \C^{\P(r)*},
\end{equation}
and so one on $\G=\Grass(S^K \C^{\P(r)*},\P(Kr))$. It is this action that
commutes with the action on $\Hilb\subset\G$ induced by that on $\PP$.
Points in the orbit of $x\in\G$ correspond to
the projective transformations of $X$.

From (\ref{eq:seshilb}), the fibre over $x\in\Hilb$ of the
hyperplane bundle on $\G$ is naturally isomorphic to
\begin{equation} \label{eq:line}
\O_{\G}(1)\res{x}=\Lambda^{\max} H^0_X(L^{Kr})\otimes
\left(\Lambda^{\max} S^K H^0_X(L^r)\right)^*.  
\end{equation}

\begin{defn}
$(X,L)$ is Hilbert (semi)stable with respect to $r$ if the point $x_{r,K}\in\Hilb$
is GIT (semi)stable for the action of $SL(\P(r),\C)$, linearised on (\ref{eq:line}),
for all $K\gg0$.
\end{defn}

\emph{Asymptotic Hilbert stability} is defined to mean Hilbert stability
for all sufficiently large $r$.
By picking a different line bundle on the Hilbert scheme (i.e.\ a different
projective embedding of Hilb -- the beautiful Chow embedding \cite{Mu})
we also get the notion of \emph{Chow stability} with respect to $r$
and \emph{asymptotic Chow stability}. We need the concept
of a \emph{test configuration}, as defined in the foundational paper \cite{Do2}.

\begin{defn}
A \emph{test configuration} for a polarised variety $(X,l)$ is
\begin{enumerate}
  \item A proper flat morphism $\pi\colon\X\to\C$,
  \item An action of $\C^\times$ on $\X$ covering the usual
    action of $\C^\times$ on $\C$,
  \item An equivariant very ample line bundle $\L$ on $\X$,
  \end{enumerate}
  such that the fibre $(\X_t,\L|_{\X_t})$ is isomorphic to $(X,l)$ for
  one, and so all, $t\in\C\backslash\{0\}$.

A test configuration is called a \emph{product configuration}
if $\X\cong X\times\C$, and a \emph{trivial configuration} if in addition
$\C^\times$ acts only on the second factor. 

We will often need a weaker concept which we call a \emph{semi test configuration}
where $\L$ is just globally generated.
\end{defn}

\begin{prop}
In the situation of (\ref{rlarge}), a 1-parameter subgroup of $GL(\P(r),\C)$
is equivalent to the data of a test configuration for $(X,L^r)$.
\end{prop}

\begin{proof}
The action of a 1-parameter subgroup of $GL(\P(r),\C)$ on
$X\subseteq\PP(H^0_X(L^r)^*)$
clearly gives a test configuration $(\X,\L)$ for $(X,L^r)$ with $\L$ the
pull back of the (equivariant) line bundle $\O_\PP(1)$.

Conversely the subgroup
can be recovered from the test configuration via the induced $\C^\times$-action
on the dual of the vector space $(\pi_*\L)|_{\{0\}}$. Since $\pi_*\L$ is
torsion-free (by flatness) over the curve $\C$ it is a vector bundle, and
$\X$ sits inside the projectivisation of its dual by the very ampleness of
$\L$. Its general fibre has dimension $\P(r)$, so $(\pi_*\L)|_{\{0\}}$ does
too.

Pick a trivialisation of $\pi_*\L$ over $\C$, identifying it with
$(\pi_*\L)|_{\{0\}}\times\C$. This has a diagonal $\C^\times$-action, inducing
one on $\PP(\pi_*\L)^*\supset\X$ which yields the original test configuration;
thus these two operations are mutual inverses.

(Note that in fact $(\pi_*\L)|_{\{0\}}\cong H^0_\X(\L)\big/tH^0_\X(\L)$.
The map $\leftarrow$ is just restriction; we need to define its inverse $\to$.
Any element of $(\pi_*\L)|_{\{0\}}$ can be lifted to a meromorphic section
$s$ of $\L$ on $\X$ that is regular on $\X_0$.
The projection of its polar locus to $\C$ does not contain $\{0\}$
and so is a finite number of points in $\C$. Choosing a polynomial $p$ with
high order zeros at these points such that $p(0)=1$, $p.s$ is a holomorphic
section with the same class as $s$ in $(\pi_*\L)|_{\{0\}}$ since $p-1\in(t)$.
Then $[p.s]$ defines the required class in $H^0_\X(\L)\big/tH^0_\X(\L)$.)
\end{proof}

Denote the weight of the induced $\C^\times$-action on
$(\pi_*\L^K)|_{\{0\}}=H^0_\X(\L^K)\big/tH^0_\X(\L^K)$ by $w(Kr)$.
(We enumerate by $k:=Kr$ since $(\X,\L^K)$ is a test configuration for $(X,L^{Kr})$.
The confused reader may set $r=1$ temporarily.) For $K\gg0,\ 
(\pi_*\L^K)|_{\{0\}}$ is $H^0_{\X_0}(\L^K)$ and
$w(k)=w(Kr)$ is a polynomial of degree $n+1$ in $k=Kr$ by the equivariant
Riemann-Roch theorem. (It is important that we do \emph{not} modify
$w(k)$ to be this polynomial for small $k$, so for instance $w(r)$ really
is the weight on
$(\pi_*\L)|_{\{0\}}$.) To make the $\C^\times$-action special
linear on $(\pi_*\L)|_{\{0\}}$ we first pull back the family by the cover
$\C\to\C,\ t\mapsto t^{r\P(r)}$, making the action of weight $r\P(r)w(r)$.
Composing with the trivial action which scales the $\L$-fibres with weight
$-rw(r)$ scales $\Lambda^{\max}(\pi_*\L)|_{\{0\}}$ with weight
$-r\P(r)w(r)$, cancelling out the previous action. (The extra factor of $r$
is to make the formula (\ref{normalised}) nicer, and is natural if $(X,\L)$
is the $r$th twist of a test configuration for $(X,L)$.)
Since this new action acts with zero total weight on $S^K(\pi_*\L)|_{\{0\}}$,
it acts on the line $\O_\G(1)|_{\{\X_0\}}=
\Lambda^{\max}H^0_{\X_0}(\L^{Kr})\otimes
\left(\Lambda^{\max}S^K(\pi_*\L)|_{\{0\}}\right)^*$ (compare (\ref{eq:line})
which was for $(X,L^r)$ with no higher cohomology of $L^r$) with
\emph{normalised weight}
\begin{equation} \label{normalised}
\tilde{w}_{r,Kr}=\tilde{w}_{r,k}=w(k)r\P(r)-w(r)k\P(k), \qquad k:=Kr.
\end{equation}
The normalised weight is a polynomial $\sum_{i=0}^{n+1}e_i(r)k^i$
of degree $n+1$ in $k$ for $k\gg0$,
with coefficients which are also polynomial of degree $n+1$ in $r$ for $r\gg0$:
$e_i(r)=\sum_{j=0}^{n+1}e_{i,j}r^j$ for $r\gg0$. The normalisation means
that the coefficient of
$(rk)^{n+1}$ vanishes: $e_{n+1,n+1}=0$. The Hilbert-Mumford
criterion relates this weight to stability as follows.

\begin{thm} \label{thm:instability}
A polarised variety $(X,L)$ is \emph{stable} if and only if
$$
\tilde{w}_{r,k} \succ 0 \quad \forall \text{ nontrivial
test configurations for } (X,L^r),
$$
where $\succ$ and $\forall$ mean the following for the
different notions of stability:
  \begin{itemize}
  \item Hilbert stable with respect to $r$: for any nontrivial
test configuration for $(X,L^r),\ \tilde{w}_{r,k}>0$ for all $k\gg0$,
\item Asymptotically Hilbert stable: for all $r\gg0$, any nontrivial
test configuration for $(X,L^r)$ has $\tilde{w}_{r,k}>0$ for all $k\gg0$,
  \item Chow stable with respect to $r$: for any nontrivial
test configuration for $(X,L^r)$, the leading $k^{n+1}$
coefficient $e_{n+1}(r)$ of $\tilde{w}_{r,k}$ is positive: $e_{n+1}(r)>0$,
\item Asymptotically Chow stable: for all $r\gg0$, any nontrivial
test configuration for $(X,L^r)$ has $e_{n+1}(r)>0$,
\item K-stable: for all $r\gg0$, any nontrivial test configuration for $(X,L^r)$
has leading
coefficient $e_{n,n+1}$ of $e_{n+1}(r)$ (the \emph{Donaldson-Futaki invariant}
of the test configuration) positive: $e_{n,n+1}>0$.
\end{itemize}
Furthermore the result holds if we replace ``stable" with
``semistable" and strict inequalities with non strict inequalities throughout.

Finally $(X,L)$ is \emph{polystable} if it is semistable and any destabilising
test configuration (i.e.\ one which is not strictly stable) is a product
configuration.
\end{thm}

The increasing number of test configurations that have to be tested in the
second, fourth and fifth definitions
as $r\to\infty$ currently prevent us from adding K-stability to the left
of the following consequences of Theorem \ref{thm:instability}. \\

Asymptotically Chow stable $\Rightarrow$ Asymptotically Hilbert stable
$\Rightarrow$ Asymptotically Hilbert semistable $\Rightarrow$ Asymptotically
Chow semistable $\Rightarrow$ K-semistable. \\

The relation between K-stability and asymptotic Chow stability is analogous
to the relation between slope stability and Gieseker stability for vector
bundles, and a geometric criterion for asymptotic Chow stability would show
it was implied by K-stability; we only have a necessary condition (Theorem
\ref{chow}).

The fact that Chow stability is controlled by the coefficient $e_{n+1}(r)$
is due to Mumford \cite{Mu}. The definition of K-stability above is due to
Donaldson, adapting Tian's original differential-geometric definition to
allow nonnormal central fibres $\X_0$ (though what
is called properly semistable in \cite{Ti2} and stable in
\cite{Do2} is what we call K-polystable). Tian \cite{Ti1}
defines a line bundle (the ``CM polarisation") on $\Hilb$ such that K-stability
is exactly stability in the sense of the Hilbert-Mumford criterion for this
line bundle \cite{PT}. K-stability is probably not a bona
fide GIT notion: Tian's polarisation may not be ample, and
the number of test configurations increases as $r$ tends to infinity.
However it is K-polystability that is conjecturally related to
the existence of constant scalar curvature K\"ahler metrics; we apply our
methods to these in \cite{RT}.

The definition of polystability says that any destabilising test configuration
comes from a $\C^\times$-action on $X$ with the appropriate weight 0 (i.e.\
Donaldson-Futaki invariant 0 in the K-polystability case). This corresponds to an
orbit in Hilb which is closed in the semistable points but with possibly
higher dimensional stabilisers; equivalently the orbit in the total space
of the dual of the polarising line bundle over Hilb is closed. There seems
to be no universally accepted name for this; we use polystability by analogy
with bundles.

A test configuration $(\X,\L)$ for $(X,L^r)$ can be twisted to get another,
$(\X,\L^K)$, for $(X,L^{Kr})$; for $K\gg0$ this will have no higher cohomology.
Since the definition of K-stability is unchanged if
$L$ is replaced by some power, for this notion we can allow $\L$ to
be an ample $\Q$\,-line bundle in the definition of a test configuration.

Letting $F:=e_{n,n+1}$ denote the Donaldson-Futaki invariant and writing
the unnormalised weight $w(k)$ as $b_0k^{n+1}+b_1k^n+\ldots\,$, we see that
\begin{equation} \label{Futaki}
F=b_0a_1-b_1a_0,
\end{equation}
and $-a_0^{-2}F$ is the coefficient of $k^{-1}$ in the expansion of
$w(k)/k\P(k)$.  When the central fibre $\X_0$ is smooth,
$F=\frac{a_0}{4}\nu$, where $\nu$ is the usual Calabi-Futaki invariant
of $c_1(L)$ and the vector field generated by the $S^1$-action
\cite{Do2}. \\

Without loss of generality we now take $r=1$.
An arbitrary test configuration $(\Y,\O_\Y(1))$ for $(X,L)$ is, by definition,
$\C^\times$-isomorphic to the trivial test configuration $(X\times\C,L)$
away from the central fibre. It is therefore $\C^\times$-birational to
$(X\times\C,L)$ and so is dominated by a blow up $(\X,\L)$ of $X\times\C$
in a $\C^\times$-invariant ideal $I$ supported on (a thickening
of) the central fibre $X\times\{0\}$:
\begin{eqnarray}
(\X,\L)\ =&(\B_I(X\times\C),L(-E))&\Rt{\phi}\ (\Y,\O_\Y(1)) \nonumber \\
&\downarrow p \label{any} \\ &(X\times\C,L)\,. \nonumber
\end{eqnarray}
Here we use the canonical $\C^\times$-action
on $\B_I(X\times\C)$ inherited from that on $I$, and its linearisation on
the line bundle $\L:=L-E=p^*(L\otimes I)$, where $E$ denotes the exceptional
divisor. $L-E=\phi^*\O_\Y(1)$ and the horizontal arrow in (\ref{any}) is
an equivariant map of equivariant
polarisations (although $L-E$ may be only semi-ample), whereas
$p$ does not respect the polarisation.

Mumford (\cite{Mu} section 3 of the proof of Theorem 2.9) essentially shows
that any test configuration is of this form, where $I=I_r$ is of the form
\begin{equation} \label{ideals}
I_r=\I_0+t\I_1+t^2\I_2+\ldots+t^{r-1}\I_{r-1}+(t^r)\ \subset\ \O_X\otimes\C[t].
\end{equation}
Here the ideals $\I_0\subseteq\I_1\subseteq\ldots\subseteq\I_r
\subseteq\O_X$ correspond to subschemes $Z_0\supseteq Z_1\supseteq\ldots\supseteq
Z_{r-1}$ of $(X,L)$, so $I_r$ is $\C^\times$-invariant and corresponds to
a subscheme of $X\times\C$ supported on (the $r$-fold
thickening of) the central fibre $X\times\{0\}$.

This can be proved by writing the $\C^\times$-invariant ideal $I$ as a sum
of weight spaces: defining $\I_j$ in terms of the weight-$j$ piece $t^j\I_j$
we get the weight space decomposition (\ref{ideals}). Or, embedding a test
configuration in some $\PP^N\times\C$, Mumford's result applies directly.

For example, given a $\C^\times$-action on $X$ with ``repulsive fixed point
set" $Z$ (that part of the fixed point set with negative weight spaces in
its normal cone), there is a blow up of $X\times\C$ supported on $Z\times\{0\}$ in which the proper transform of $X\times\{0\}$ can be blown down to give
back $X\times\C$ but with a nontrivial $\C^\times$-action.

Since $(\X,\L)$ is a \emph{semi} test configuration (and
because later we will want to replace an ample $(X,L)$ with a semi-ample
$(\widehat X,p^*L)$ for some resolution of singularities $p$) we will
prove many results in the generality of semi test configurations.

General test configurations (\ref{any}) will be studied in Section \ref{testconfigs};
next we look at the simplest case (beyond product configurations) of $r=1,\
I=\I_Z+(t)$ in (\ref{any}): the \emph{deformation of $X$ to the normal cone
of $Z$}.

\section{Deformation to the normal cone and K-slope stability} \label{Kslope}

Given any proper subscheme $Z\subset X$ we get a canonical test
configuration $\X$, the blow up of $X\times\C$ along $Z\times\{0\}$ with
exceptional divisor $P$.  This deformation to the normal cone has central
fibre $\X_0=\widehat X\cup_E P$, where $\widehat X$ is
the blow up of $X$ along $Z$ with exceptional divisor $E$. When $Z$ and
$X$ are both smooth $P=\PP(\nu\oplus\underline\C)$ is the projective
completion of the normal bundle $\nu$ of $Z$ in $X$, glued along
$E=\PP(\nu)$ to the blow up of $X$.  For more on the deformation to the
normal cone see \cite{Fu}; for a diagram see Section \ref{testconfigs}.
 
Let $\pi$ be the composite of the projections
$\X\to X\times\C\to X$. For $L$ an ample line bundle on $X$ and $c$ a positive
rational number let $\L_c$ be the $\Q\,$-line bundle $\pi^*L-cP$. This
restricts to $L$ on the general fibre of $\X\rightarrow\C$, and is ample
for $c$ sufficiently small.

\begin{prop} \label{ample}
Fix $L$ an ample (respectively semi-ample) line bundle on $X$.
  For $c\in(0,\epsilon(Z))\cap\Q$ (\ref{def:seshadri}), the $\Q$-line bundle
  $\L_c$ is ample (semi-ample) on
  $\X$. If $\epsilon(Z)\in\Q$ then $\L_{\epsilon(Z)}$ is nef. If in
  addition the global sections of $L^k\otimes\I_{\!Z}^{\epsilon(Z)k}$
  saturate (see Section \ref{notation}) for $k\gg0$ then $\L_{\epsilon(Z)}$
  is semi-ample.
\end{prop}

\begin{proof}
Note that if $k,kc\in\mathbb N$ and $L^k$ is globally generated
and the sections of $L^k\otimes\I_{\!Z}^{ck}$ saturate, then $\L_c^k$
is globally generated on $X\times\C$ by the
sections $\pi^*H^0(L^k\otimes\I_{\!Z}^{ck})+t^{ck}\pi^*H^0(L^k)$. Algebraically
this is the statement that on $X\times\C$, $\I_{\!Z}^{ck}+(t^{ck})$ saturates
$(\I_Z+(t))^{ck}$. Geometrically it says that the global sections saturating
$L^k\otimes\I_{\!Z}^{ck}$ generate $L^k-ckP$ on $\X$ away from the zero
section $\PP(\underline\C\to Z)\cong Z$ of the cone $P\to Z$, over which
the sections $t^{ck}H^0(L^k)$ generate.

Putting $c=\epsilon$ now gives the third claim. For the first
two we claim that we may assume that $L$ is ample by replacing $(X,L)$ by
$(Y,\O_Y(1)):=\Proj\bigoplus_kH^0(L^k)$ if necessary. For $k$ sufficiently
large that $L^k\otimes\I_{\!Z}$ is globally
generated, $H^0(L^k\otimes\I_{\!Z})\subset H^0(L^k)=H^0(\O_Y(k))$
generates a subsheaf $\O_Y(k)\otimes\I_{Z_0}\subset\O_Y(k)$ and so a subscheme
$Z_0\subset Y$ such the pullback of $\I_{Z_0}$ to $X$ is $\I_Z$. Thus it
is sufficient to prove the result for $Z_0\subset(Y,\O_Y(1))$ and then pullback
to get the result for $Z\subset(X,L)$.

On $\widehat X$, $L-cE$ is in the ample cone for small $c$ and
on its boundary for $c=\epsilon(Z)$, so by Kleiman \cite{Kl}
$L-cE$ is ample for rational $c<\epsilon(Z)$. So $L^k-ckE$ is globally generated
for $k\gg0$; equivalently the global sections of $L^k\otimes\I_{\!Z}^{ck}$
saturate $\I_{\!Z}^{ck}$ (for $k$ sufficiently large that the pushdown of
$\O(-ckE)$ to $X$ is $\I_{\!Z}^{ck}$). $L^k$ is also globally generated so
again this implies that $\L_c^k$ is globally generated. Thus $\L_c$
is nef for $c\in(0,\epsilon(Z))$, but it is ample for small $c$ since $\pi^*L$
is ample. Thus by Kleiman again, $\L_c$ is actually ample and $\L_{\epsilon(Z)}$
is nef.
\end{proof}

The obvious $\C^\times$-action on $X\times\C$ (acting trivially on the
$X$ factor) lifts to an action on the deformation to the normal cone
$(\X,\L_c)$ which, on the central fibre $\X_0=\widehat
X\cup_EP$, is trivial on $\widehat X$.  When $Z$ and $X$
are both smooth, $\lambda\in\C^\times$ acts on
$P=\PP(\nu\oplus\underline \C)$ as $\operatorname{diag}(1,\lambda)$.

\begin{thm} \label{degen}
Fix $L$ ample and $c\in(0,\epsilon(Z))\cap\Q$. Then $(\X,\L_c)$ is a flat
family of polarised schemes, and, for all $k\gg0,\ ck\in\mathbb N$, the
total weight of the induced action on $H^0(\X_0,\L_c^k)$ is
$$
w(k)=-\sum_{j=1}^{ck} j\,h^0\!\left(L^k\otimes(\I_{\!Z}^{ck-j}/
\I_{\!Z}^{ck-j+1})\right)=\sum_{j=1}^{ck} \chi(L^k\otimes
\I_{\!Z}^j)-ckh^0(L^k).
$$
\end{thm}

\begin{proof} From the definition of the blow up of a subscheme, $\X=\Proj
\bigoplus_{k\ge 0} \S_k$, where
\begin{eqnarray} \nonumber
\S_k &=& (\I_{\!Z \times \{0\}})^k = (\I_Z+(t))^k \\
&=& \bigoplus_{j=0}^{k-1} t^j\I_{\!Z}^{k-j}\ \oplus\
t^k\C[t]\O_X\ \subset\ \C[t]\otimes\O_X. \label{S}
\end{eqnarray}
It follows that for all $k\gg0$, the pushdown of $-kP$ to $X\times\C$ is
$\bigoplus_{j=0}^{k-1}t^j\I_{\!Z}^{k-j}\oplus t^k\C[t]\O_X$, with the higher
pushdowns zero. By the ampleness of $L-cP$ (\ref{ample}), for $k\gg0$ we
have the vanishing of
$$
H^i_\X((L-cP)^k) = H^i_{X\times\C}(\pi_*(L^k-ckP)) =
\bigoplus_{j=0}^{ck-1} t^j H^i_X(L^k\otimes\I_{\!Z}^{ck-j})\oplus
t^{ck}\C[t]H^i_X(L^k),
$$
for $i>0$. In particular then,
\begin{equation} \label{vanish}
H^i_X(L^k\otimes\I_{\!Z}^{ck-j})=0\quad\text{for}\ j=0,\ldots,ck, \quad
i>0.
\end{equation}

Now
$$
\X_0=\Proj\bigoplus_{k\ge0}\S_k/t\S_k,
$$
where by (\ref{S}),
\begin{equation} \label{St}
\S_k/t\S_k=\I_{\!Z}^k\ \oplus\ t\!\left(\!\I_{\!Z}^{k-1}/\I_{\!Z}^k\right)\
\oplus\,\ldots\,\oplus\ t^k\big(\O_X/\I_Z\big).
\end{equation}
Replacing $\I_Z$ by $\I^c_{\!Z}$ does not change the blow up (just the corresponding
exceptional line bundle) so that, for $k\gg0$,
\begin{equation} \label{wspace}
H^0_{\X_0}(\L_c^k)=H^0_X(L^k\otimes\I_{\!Z}^{ck})\ \oplus\,\bigoplus_{j=1}^{ck}
t^j H^0_X(L^k\otimes(\I_{\!Z}^{ck-j}/\I_{\!Z}^{ck-j+1})).
\end{equation}
The vanishing (\ref{vanish}) of
$H^1(L^k\otimes\I_{\!Z}^{ck-j+1})$ for $1\le j\le ck$ means that the dimension
of this is
$$
h^0_{\X_0}(\L_c^k)=
h^0_X(L^k\otimes\I_{\!Z}^{ck})+\sum_{j=1}^{ck}\left(h^0_X(L^k\otimes \I_{\!Z}^{ck-j})-
h^0_X(L^k\otimes\I_{\!Z}^{ck-j+1})\right),
$$
which is $h^0_X(L^k)$. By (\cite{Ha} Theorem III.9.9) this proves flatness
of the family.

Now (\ref{wspace}) is the weight space decomposition with respect to the
$\C^\times$-action: $\C^\times$ acts on $t$ with weight $-1$ and trivially
on $X$ and so on $\I_Z$. Therefore
\begin{eqnarray}
w(k) &=& -\sum_{j=1}^{ck} jh^0\left(L^k\otimes(\I_{\!Z}^{ck-j}/
\I_{\!Z}^{ck-j+1})\right) \nonumber \\
&=&-\sum_{j=1}^{ck} j\left(h^0(L^k\otimes \I_{\!Z}^{ck-j})-
h^0(L^k\otimes \I_{\!Z}^{ck-j+1})\right) \label{noco} \\
&=& \sum_{j=1}^{ck} (ck-j+1) h^0(L^k\otimes\I_{\!Z}^j)-\sum_{j=0}^{ck-1}
(ck-j)h^0(L^k\otimes \I_{\!Z}^j) \nonumber \\
&=& \sum_{j=1}^{ck} h^0(L^k\otimes \I_{\!Z}^j)-ckh^0(L^k)
\ =\ \sum_{j=1}^{ck} \chi(L^k\otimes \I_{\!Z}^j)-ckh^0(L^k), \nonumber
\end{eqnarray}
where the second and last equalities follow from the vanishing (\ref{vanish})
of $H^i(L^k\otimes\I_{\!Z}^{ck-j+1})$ for $1\le j\le ck$ and $i>0$.
\end{proof}

We may only take $c=\epsilon(Z)$ if the global sections
of $L^k\otimes\I_{\!Z}^{\epsilon(Z)k}$ saturate for $k\gg0$. In this case
$\L_c$ is semi-ample, pulled back from a contraction $p$ from $(\X,\L_c)\to\C$
to $\Proj\bigoplus_{k\gg0}H^0_\X(\L_c^k)$, which we call $(\Y,\O_\Y(1))\to\C$.
By construction, $p^*\colon H^0_\Y(\O_\Y(k))\to H^0_\X(\L_c^k)$ is then
an isomorphism. By Lemma 2.13 of \cite{Mu}, $\Y\to\C$ is also a flat family,
so $(\Y,\O_\Y(k))\to\C$ is a test configuration for $k\gg0$.

\begin{thm} \label{ss}
Let $Z$ be a proper subscheme of an irreducible polarised algebraic
variety $(X,L)$, and suppose that $c=\epsilon(Z)\in\Q$ and the global sections
of $L^k\otimes\I_{\!Z}^{\epsilon(Z)k}$ saturate for $k\gg0,\ ck\in\mathbb
N$. Letting $(\Y,\O_\Y(1))$ be the contraction of $(\X,\L_c)$ constructed
above, the weight of the action on $\Lambda^{\max}H^0_{\Y_0}(\O_\Y(k))$ is
$$
w(k)=\sum_{j=1}^{ck}\chi(L^k\otimes\I_{\!Z}^j)-ckh^0(L^k)=
\sum_{j=1}^{ck}\chi(L^k\otimes\I_{\!Z}^j)-ckh^0(L^k)+O(k^{n-1}).
$$
\end{thm}

\begin{proof}
By (\ref{S}), for $k\gg0$,
$$
H^0_\X(\L_c^k)=H^0_{X\times\C}(L^k\otimes\S_{ck})\cong
\bigoplus_{j=0}^{ck-1}t^jH^0_X(L^k\otimes\I_{\!Z}^{k-j})\ \oplus\ t^{ck}\C[t]H^0_X(L^k),
$$
yielding
\begin{equation} \label{qt}
\frac{H^0_\X(\L_c^k)}{tH^0_\X(\L_c^k)}\cong H^0_X(L^k\otimes\I_{\!Z}^{ck})\
\oplus\ \bigoplus_{j=1}^{ck}\,t^j\frac{H^0_X(L^k\otimes\I_{\!Z}^{ck-j})}
{H^0_X(L^k\otimes\I_{\!Z}^{ck-j+1})}\,.
\end{equation}
By the isomorphism $p^*\colon H^0_\Y(\O_\Y(k))\to H^0_\X(\L_c^k)$ this is
the weight space decomposition of $H^0_\Y(\O_\Y(k))\big/tH^0_\Y(\O_\Y(k))$,
which by flatness of $\Y\to\C$ and ampleness of $\O_\Y(1)$ is
$H^0_{\Y_0}(\O_\Y(k))$. So the total weight is
$$
-\sum_{j=1}^{ck} j\left(h^0(L^k\otimes \I_{\!Z}^{ck-j})-
h^0(L^k\otimes \I_{\!Z}^{ck-j+1})\right).
$$
Just as in (\ref{noco}) this is $\sum_{j=1}^{ck}h^0(L^k\otimes\I_{\!Z}^j)-ckh^0(L^k)$.
Then replacing $h^0$ by $\chi$ gives an error bounded by $\sum_{j=1}^{ck}
h^{\ge1}_X(L^k\otimes\I_{\!Z}^{ck-j+1})$, where $h^{\ge1}:=\sum_{i=1}^nh^i$.
So the result follows from Lemma \ref{error} below.
\end{proof}

\begin{lem} \label{error}
Fix $Z\subset(X,L)$ and $c\in(0,\epsilon(Z))\cap\Q$, or $c\in(0,\epsilon(Z)]
\cap\Q$ if $\epsilon(Z)\in\Q$ and the global sections of $L^k\otimes
\I_{\!Z}^{\epsilon(Z)k}$ saturate for $k\gg0$. Then
$$
\sum_{j=0}^{ck}h^{\ge1}_X(L^k\otimes\I_{\!Z}^j)=O(k^{n-1}).
$$
\end{lem}

\begin{proof}
Fix $\varepsilon\in(0,c)\cap\Q$. Then $\L_{c-\varepsilon}$ is ample on
$\X$, while $\L_c$ is generated by its global sections and so nef. So we
can apply Fujita vanishing (\cite{La} Theorem
1.4.35) to give $N\gg0$ such that $\L_{c-\varepsilon}^N\otimes\L_c^p$ has
no higher cohomology for any $p\ge0$. For $k>N$, setting $p=k-N$ shows
that $\L_{c-\varepsilon N/k}^k$ has no higher cohomology. So for $k\gg0$
we have
\begin{eqnarray} \nonumber
0&=&H^i_\X((L-(c-\varepsilon N/k)P)^k)=H^i_{X\times\C}
(\pi_*(L^k-(ck-\varepsilon N)P)) \\
&=&\bigoplus_{j=0}^{ck-\varepsilon N-1}t^j H^i_X(L^k\otimes
\I_{\!Z}^{ck-\varepsilon N-j})\ 
\oplus\ t^{ck-\varepsilon N}\C[t]H^i_X(L^k), \label{dec}
\end{eqnarray}
for $i>0$. That is, $h^{\ge1}(L^k\otimes\I_{\!Z}^{ck-\varepsilon N-j})=0$ for
$j=0,\ldots,ck-\varepsilon N$. So the sum becomes
$$
\sum_{j=0}^{ck}h^{\ge1}_X(L^k\otimes\I_{\!Z}^j)=
\sum_{j=0}^{\varepsilon N-1}h^{\ge1}_X(L^k\otimes\I_{\!Z}^{ck-j})=
\sum_{j=0}^{\varepsilon N-1}h^{\ge1}_{\widehat X}(L^k-(ck-j)E),
$$
where we have taken $k$ sufficiently large that the pushdown of $\O_{\widehat
X}(-iE)$ is $\I_{\!Z}^i$ (and higher pushdowns are zero) for $i>ck-\varepsilon
N$. A corollary of Fujita vanishing is that $h^i(\curly F(kD))=O(k^{n-i})$
for any coherent sheaf $\curly F$ and nef divisor $D$ (\cite{La} Theorem
1.4.40). Applying this on $\widehat X$ in turn to $\curly F=
\O,\O(E),\ldots,\O((\varepsilon N-1)E)$ and $D=L-cE$ shows that each
$h^{\ge1}_{\widehat X}(L^k-(ck-j+1)E)=O(k^{n-1})$. Since
$N$ is fixed, then, we get a similar bound on the whole sum.
\end{proof}

\begin{cor} \label{weight}
For $c\in(0,\epsilon(Z)]\cap\Q$ define
$$
\tilde{w}_{r,k}(c)=r\chi(L^r)\sum_{j=1}^{ck}h^0(L^k\otimes\I_{\!Z}^j)
\,-k\chi(L^k)w(r)-ck\chi(L^k)r\chi(L^r),
$$
which, for $r$ sufficiently large \emph{for fixed $Z\subset X$,} is
$$
\tilde{w}_{r,k}(c)=r\chi(L^r)\sum_{j=1}^{ck}h^0(L^k\otimes
\I_{\!Z}^j) - k\chi(L^k)\sum_{j=1}^{cr}h^0(L^r\otimes\I_{\!Z}^j).
$$
For $cr\in\mathbb Z$, $L^r$ globally generated and $L^r\otimes\I_{\!Z}^{cr}$
saturated by its sections, $(X,L^r)$ is unstable if $\tilde{w}_{r,k}(c)\preceq0$,
where $\preceq$ is defined as in Theorem \ref{thm:instability}, depending
on the type of stability. We say that $(X,L^r)$ is destabilised by $Z$. Similarly
for strict instability and $\prec$.
\end{cor}

\begin{proof}
The conditions on $r$ are just to ensure that $(\X,\L_c^r)$ (or $(\Y,\O_\Y(r))$
if $c=\epsilon$) is a genuine test configuration; for K-stability we are
free to twist by higher $r$ and use Proposition \ref{ample}
to remove these conditions.

The result is a direct consequence of Theorems \ref{degen} and \ref{ss} and
the identity $\tilde{w}_{r,k}=w(k)r\chi(L^r)-w(r)k\chi(L^k)$ (\ref{normalised}).
\end{proof}

Recall from the introduction the definition of the $a_i$ and $a_i(x)$ (\ref{HS}).
Choose $p$ so that $\mathbf R\pi_*\O(-jE)=\I_{\!Z}^j$ for $j\ge p$,
then for $k,xk\in\mathbb N,\,k\ge p/x$,
\begin{equation} \label{definitionofa_i}
\chi\_{\widehat X}(L^k-xkE)=\chi\_X(L^k\otimes\I_{\!Z}^{xk})=a_0(x)k^n
+ a_1(x)k^{n-1} + \ldots +a_n(x).
\end{equation}
But by Riemann-Roch, $P(k,j):=\chi\_{\widehat{X}}(L^k-jE)$ is a polynomial
of total degree $n$. Writing $P=P_0+\ldots+P_n$ where $P_i$ is
homogeneous of degree $n-i,\ a_i(x)=P_i(1,x)$ is therefore a degree
$n-i$ polynomial in $x$ which can be extended to
all real $x$. 

\begin{prop} \label{integrals}
For $X$ irreducible the weights $\sum_{j=1}^{ck}\chi(L^k\otimes \I_{\!Z}^j)-
ckh^0(L^k)$ of Theorems \ref{ss} and \ref{degen} can be written
$$
w(k)=\left(\int_0^ca_0(x)dx-ca_0\right)k^{n+1}+\left(\int_0^c\!a_1(x)
+\frac{a_0'(x)}2\,dx-ca_1\right)k^n+O(k^{n-1}).
$$
\end{prop}

\begin{proof} Using (\ref{definitionofa_i}) for $j>p$ we split up
$\sum_{j=1}^{ck}\chi(L^k\otimes \I_{\!Z}^j)$ as
\begin{eqnarray} \nonumber
&& \sum_{j=1}^p\chi(L^k\otimes \I_{\!Z}^j)
+\sum_{j=p+1}^{ck}\chi(L^k\otimes \I_{\!Z}^j) \\
&=&p\,\chi(L^k)-\sum_{j=1}^p\chi(L^k\otimes\O_{jZ})+\sum_{j=p+1}^{ck}(a_0(j/k)k^n
+ a_1(j/k)k^{n-1}+\ldots) \nonumber \\
&=&\int_0^{p/k}\!a_0(x)dx\,k^{n+1} + \sum_{j=p+1}^{ck}(a_0(j/k)k^n
+ a_1(j/k)k^{n-1})+O(k^{n-1}), \label{prox}
\end{eqnarray}
using the fact that $jZ$ has dimension $\le n-1$ (since $X$ is irreducible)
and $a_0(x)$ is a polynomial with $a_0(0)=a_0$.

For a smooth function $f$, as $k\to\infty$ with $ck\in\mathbb
N$, the trapezium rule gives \cite{Hi}
\begin{equation} \label{easyeuler}
\sum_{j=1}^{ck}f(j/k)=\int_0^c\left(kf(x)+\frac{f'(x)}2\right)dx+O(k^{-1}).
\end{equation}
This can be proved by Taylor's theorem, or directly for polynomials by noting
that for $f(x)=x^m$ we have the identity 
$$
\sum_{j=1}^{ck} j^m = \frac{(ck)^{m+1}}{m+1} + \frac{(ck)^m}2+ O(k^{m-1}).
$$
Applying this to $f(x)=a_i(x),\ i=1,2$, (\ref{prox}) now approximates
$\sum_{j=1}^{ck}\chi(L^k\otimes \I_{\!Z}^j)$ by
\begin{multline*} \hspace{-4mm}
\int_0^{p/k}\!\!a_0(x)dx\,k^{n+1}
+\left(\int_{p/k}^cka_0(x)+\frac{a_0'(x)}2dx\right)k^n+
\left(\int_{p/k}^cka_1(x)dx\right)k^{n-1}+ O(k^{n-1}) \\
=\left(\int_0^c a_0(x) dx \right) k^{n+1} + \left(\int_0^c
  a_1(x) + \frac{a_0'(x)}{2}dx \right) k^n + O(k^{n-1}),
\end{multline*}
since $\int_0^{p/k}\!a_1(x)+\frac{a_0'(x)}{2} dx\,k^n=O(k^{n-1})$.
Subtracting $ckh^0(L^k)=ca_0\,k^{n+1}+ca_1\,k^n+O(k^{n-1})$ gives the result.
\end{proof}

Note that by Riemann-Roch on $\widehat X$,
$n!a_0(x)=(L-xE)^n>0$ for $x\in(0,\epsilon(Z))$ by the ampleness of $L-xE$,
so for $c\in(0,\epsilon(Z)]$, $\int_0^ca_0(x)dx>0$. Therefore we can define
the \emph{slope} (or K-slope) of $\I_Z$ (\ref{def:slope}) by 
$$
\mu_c(\I_Z)=\mu_c(\I_Z,L) = \frac{\int_0^c a_1(x) + \frac{a'_0
    (x)}{2} dx} {\int_0^c a_0(x)
  dx}=\frac{\int_0^ca_1(x)dx}{\int_0^ca_0(x)dx}+
\frac{a_0(c)-a_0}{2\int_0^ca_0(x)dx}\,,
$$
and that of $X$ (\ref{def:X}),
$$
\mu(X) = \frac{a_1}{a_0} \qquad \left(=-\frac n2\cdot\frac{K_X.L^{n-1}}{L^n}
\ \text{ for $X$ smooth}\right).
$$
This gives the following obvious definition of K-slope stability, which,
like K-stability, is independent of replacing
$L$ by $L^r$ (on replacing $c,\,x$ and $\epsilon$ by $rc,\,rx$ and
$r\epsilon$).

\begin{defn} \label{defstab}
We say that $(X,L)$ is slope (semi/poly)stable if for every proper
subscheme $Z$ of $X$, the following holds:
\begin{itemize}
\item \textbf{\emph{slope semistability}}: $\mu_c(\I_Z,L)\le\mu(X)$ for all
$c\in(0,\epsilon(Z)]$,
\item \textbf{\emph{slope stability}}: \hspace{6mm} $\mu_c(\I_Z,L)<\mu(X)$
for all $c\in(0,\epsilon(Z))$, and also for $c=\epsilon(Z)$ if
$\epsilon(Z)\in\Q$ and global sections of $L^k\otimes\I_{\!Z}^{\epsilon(Z)k}$
saturate for $k\gg0$,
\item \textbf{\emph{slope polystability}}: $\ (X,L)$ is slope
  semistable, and if $(Z,c)$ is any pair such that
  $\mu(\I_Z,c)=\mu(X)$, then $c=\epsilon(Z)$ is rational and, on the
  deformation to the normal cone of $Z$, $\L_c$ is pulled back from a
  product test configuration $(\Y,\O_\Y(1))\cong(X\times\C,L)$.
\end{itemize}
\end{defn}

\begin{thm}
\label{thm:kstableslopestable}
If a polarised variety $(X,L)$ is K-stable then it is slope stable.
If it is K-polystable (respectively K-semistable) then it is slope polystable
(semistable).
\end{thm}

\begin{proof}
From the deformation to the normal cone of $Z$ we get test configurations
$(\X,\L_c^k)$ (for $c<\epsilon(Z),\ k\gg0$), or $(\Y,\O_\Y(k))$ (for $c=\epsilon(Z)$
if $\epsilon(Z)\in\Q$ and the global sections of $L^k\otimes\I_{\!Z}^{\epsilon(Z)k}$
saturate). Its total weight $w(k)$ is given by Theorems \ref{degen}
and \ref{ss} respectively. Writing $w(k)=b_0k^{n+1}+b_1k^n + O(k^{n-1})$,
its Donaldson-Futaki invariant (\ref{Futaki}) is, by Proposition \ref{integrals},
\begin{eqnarray} \nonumber
F(Z) &=& b_0a_1-b_1a_0=\left(\int_0^c\!a_0(x)dx-ca_0\right)\!a_1-\left(
\int_0^c\!a_1(x)+\frac{a_0'(x)}2\,dx-ca_1\right)\!a_0 \\
&=& a_0\Big(\mu(X)-\mu_c(\I_Z)\Big)\int_0^c\!a_0(x)dx. \label{F1}
\end{eqnarray}
This is a strictly positive multiple of $\mu(X)-\mu_c(\I_Z)$, and K-stability
(K-semistability) implies it is strictly positive (nonnegative). This gives
the result so long as, in the semistable/polystable case, the test configuration
is not trivial. But the central fibres of both $(\X,\L_c)$ and $(\Y,\O_\Y(1))$
have nontrivial $\C^\times$-actions.
\end{proof}

This result means that slope instability provides an obstruction to the existence
of constant scalar curvature K\"ahler metrics (and so also K\"ahler-Einstein
metrics); see \cite{RT} where there are also numerous examples.

Letting $\tilde{a}_i(x) = a_i-a_i(x)$ be the coefficients of
$\chi(L^k\otimes\O_{xkZ})=\chi(L^k/(L^k\otimes\I^{xk}_{\!Z}))$, we define
the slope of $\O_Z$ (in slightly misleading notation) to be
\begin{equation} \label{qslope}
\mu_c(\O_Z)=\frac{\int_0^c \tilde{a}_1(x) + \frac{\tilde{a}'_0 (x)}{2}
dx}{\int_0^c \tilde{a}_0(x)dx}\,.
\end{equation}
Notice that we can rephrase slope stability in the equivalent ways
$$
\mu_c(\I_Z)<\mu(X)\ \Longleftrightarrow\ \mu(X)<\mu_c(\O_Z)\ \Longleftrightarrow\
\mu_c(\I_Z)<\mu_c(\O_Z),
$$
due to the implications
$$
\frac AB< \frac CD\ \Longleftrightarrow\ 
\frac CD<\frac{C-A}{D-B}\ \Longleftrightarrow\ \frac AB<\frac{C-A}{D-B}
$$
for $0<B<D$, on setting $B=\int_0^c a_0(x) dx$ and $D=ca_0$ (so $D-B=\int_0^c
\tilde{a}_0(x)dx>0$).

\begin{rmks} \label{normal}
On the blow up $\widehat X$ we have the formula
\begin{equation} \label{a0}
a_0(x)=\frac1{n!}(L-xE)^n,
\end{equation}
and so $a_0(0)=a_0$, since the intersection $L^n$ can be calculated
equally on $X$ or $\widehat X$.

 The $a_i(0)$ are the coefficients of the polynomial
  $$
  \chi\_{\widehat{X}}(L^k)=\chi\_X ((\mathbf R\pi_*\O_{\widehat X})\otimes
  L^k)=\sum_{i=0}^n(-1)^ih^0_X(R^i\pi_*\O_{\widehat X}\otimes L^k)
  \quad\text{for }k\gg0.
  $$
  \emph{If $X$ is normal}, then $\pi_*\O_{\widehat X}=\O_X$ (\cite{Ha} proof
  of Corollary III.11.4). $E$ has dimension $n-1$ and $R^i\pi_*\O$
  is supported on points over which the fibre has dimension $\ge i$, so the
  support of $R^i\pi_* \O_{\widehat X}$ has codimension at least
  $i+1$. Hence $\chi(R^i\pi_* \O_{\widehat X}\otimes L^k) =
  O(k^{n-1-i})$, and $\chi\_{\widehat X}(\pi^*L^k)= \chi\_X(L^k) +
  O(k^{n-2})$ so $a_1(0)=a_1$ also.
  
  (The same argument shows that if $Z$ has dimension $j$, then $a_i(0)=a_i$
for $i\le\max\{n-j-1,0\}$ (and $i\le\max\{n-j-1,1\}$ for $X$ normal).
  If $Z\subset X$ are both smooth in a neighbourhood
  of $Z$, then $\mathbf R\pi_* \O_{\widehat X}=\O_X$ so $a_i(0)=a_i$ for
all $i$.)

Since $H^0(L^r\otimes\I^{xr}_{\!Z})\subseteq H^0(L^r\otimes\I^{yr}_{\!Z})$
for $x>y$, $a_0(x)$ is a decreasing function in $x$: $a_0'(x)\le0$. In fact
from (\ref{a0}),
\begin{equation} \label{a0'}
a_0'(x)=-\frac{1}{(n-1)!}
(L-xE)^{n-1}.E<0 \quad\text{for }x\in(0,\epsilon(Z)),
\end{equation}
by the ampleness of $L-xE$.
In particular $a_0(x)<a_0$ for $x\in(0,\epsilon(Z))$, showing that $\mu_c(\O_Z)$
(\ref{qslope}) is finite.

  For $X$ normal then,
  \begin{equation} \label{decreasing}
  \mu_c(\I_Z)=\mu(X)+\frac{a_0'(0)}{2a_0} + O(c)
  \end{equation}
  is strictly less than $\mu(X)$ for small $c$, and the slope inequality
  is automatically satisfied. In all of the examples we have
considered \cite{Ro, RT}, one need only test
the slope inequality at $c=\epsilon(Z)$.
  If this held in general it would
  simplify the definition of stability. Sz\'ekelyhidi \cite{Sz} has shown
  by example that for the modification of K-stability relevant to extremal
  metrics this is not the case.
\end{rmks}

\subsection*{\textbf{Simplifying destabilising subschemes}}

\begin{prop}\label{simplify}
  Suppose that $Z\subset X$ is a strictly destabilising subscheme in the
  sense that it violates the slope inequality (\ref{defstab}). Then at least
  one of the connected
  components of $Z$ strictly destabilises. Similarly if $Z$ is a thickening
$Z=mZ'$ of $Z'\subset X$ then $Z'$ strictly destabilises.
\end{prop}

\begin{proof}
  Suppose $Z=Z_1\cup Z_2$ with $Z_1$ and $Z_2$ disjoint.  Then
  $\tilde{a}_i^{Z}(x)=\tilde{a}_i^{Z_1}(x) + \tilde{a}_i^{Z_2}(x)$
  (\ref{qslope}) since $\O_{xkZ}=\O_{xkZ_1}\oplus \O_{xkZ_2}$.
  Suppose $Z$ is strictly destabilises.  Then there is a
  $c\in(0,\epsilon(Z)]\cap\Q$ such that
\begin{equation} \label{simplify2}
\mu_c(\O_Z) = \frac{\big(\int_0^c\tilde a_0^{Z_1}(x)dx\big)\mu_c(\O_{Z_1})+
\big(\int_0^c\tilde a_0^{Z_2}(x)dx\big)\mu_c(\O_{Z_2})}{\int_0^c
\tilde{a}_0^{Z_1}(x) dx + \int_0^c \tilde{a}_0^{Z_2}(x) dx} < \mu(X).
\end{equation}
This implies that for some $j\in\{0,1\}$, $\mu_c(\O_{Z_j})< \mu(X)$,
and by Lemma \ref{seshadridisconnected} $c\le \epsilon(Z_j)$, so $Z_j$ is
strictly destabilising.

Finally if $(\I_{\!Z}^m,c)$ destabilises then so does $(\I_{\!Z},mc)$ since
$\epsilon(mZ)=\frac1m\epsilon(Z)$ and
\begin{equation*}
\mu_c(\I_{\!Z}^m)=\mu_{mc}(\I_Z) + (m-1)\frac{\int_0^{mc}a_0'(x)dx}
{2\int_0^{mc}a_0(x)dx}<\mu_{mc}(\I_Z),
\end{equation*}
as $a_0'(x)<0$ for $x\in(0,mc)$ (\ref{a0'}).
\end{proof}

\begin{lem}\label{seshadridisconnected}
If $Z_1\cap Z_2=\emptyset$ then $\epsilon(Z_1\cup Z_2)\le
\min(\epsilon(Z_1),\epsilon(Z_2))$.
\end{lem}
 
\begin{proof} 
Let $\pi\colon \widehat X\to X$ be the blowup of $X$ along $Z_1\cup
Z_2$ with exceptional divisor $E=E_1\cup E_2$, where $E_i$ is the
subset of $E$ sitting over $Z_i$.  Let $\epsilon=\epsilon(Z_1\cup
Z_2)$, so by definition $\pi^* L-\epsilon E$ is nef.  If $C$ is an
irreducible curve contained in $E_2$ then $(\pi^* L-\epsilon
E_1).C\ge0$ by ampleness of $L$ and the fact that $E_1.C=0$.  If $C$
is not contained in $E_2$ then $(\pi^* L-\epsilon E_1).C=(\pi^*
L-\epsilon E).C + \epsilon E_2.C \ge 0$.  Hence by the Kleiman
criterion, $\pi^* L-\epsilon E_1$ is nef. But this line bundle is the pullback
of $L-\epsilon E_1$ from $\B_{Z_1}X$, so the latter line bundle is also nef
(\cite{La} Example 1.4.4(ii)), proving that $\epsilon\le\epsilon(Z_1)$.
\end{proof}

\begin{lem} \label{pts}
Let $(X,L)$ be a smooth polarised variety of dimension $n$ and
$\epsilon=\epsilon(p,L)$ be the Seshadri constant of some point $p$
in $X$. Then $p$ strictly destabilises if and only if 
$$
\big((-K_X).L^{n-1}\big)\epsilon(p,L)>(n+1)L^n.
$$
\end{lem}

\begin{proof}
We have $\tilde{a}_0(x)=a_0-a_0(x)=\frac{1}{n!}(L^n-(L-xE)^n)=
-\frac{(-x)^n}{n!}E^n=\frac{x^n}{n!}$ since $c_1(L)^{\cap j},\,j>0$ can be
represented by a cycle on $X$ missing $p$, so by a cycle on $\widehat X$
missing $E$. Similarly using $K_{\widehat X}=K_X+(n-1)E$,
$$
\tilde{a}_1(x)=a_1-a_1(x)=-\frac{K_XL^{n-1}-(K_X+(n-1)E)(L-xE)^{n-1}}{2(n-1)!}=
\frac{(n-1)x^{n-1}}{2(n-1)!}\,,
$$
yielding (\ref{qslope}),
$$
\mu_c(\O_p,L) = \frac{\int_0^c\frac{(n-1)x^{n-1}}{2(n-1)!}+\frac{x^{n-1}}{2(n-1)!}dx}
{\int_0^c\frac{x^n}{n!}dx}=\frac{n(n+1)}{2c}\,.
$$
So $p$ strictly destabilises if and only if $n(n+1)a_0<2ca_1$. As
this is linear in $c$ it holds for some $c\le\epsilon$ if and only if
it holds for $c=\epsilon$. Substituting $a_0=\frac1{n!}L^n$ and $2a_1=
-\frac1{(n-1)!}K_X.L^{n-1}$ gives the result.
\end{proof}

One can of course also calculate the slope of a smooth point more directly
by working locally with
$H^0(L^k\otimes\O_{ckp})\cong\bigoplus_{j=0}^{ck-1}S^jT^*_pX$. We now get

\begin{thm} \label{pts2}
If $X$ is smooth then no point strictly destabilises.
\end{thm}

\begin{proof}
For any line bundle $A$ we say that $H^0(A)$ \emph{generates $s$-jets at
$p$} if
$$
H^0(A) \rightarrow A \otimes\big(\O_p/\I_p^{s+1}\big)=A|_{(s+1)\{p\}}
$$
is a surjection, where $\I_p$ denotes the maximal ideal at $p$. We define
$s(A)$ to be the maximum integer $s$ such that $H^0(A)$ generates $s$-jets
at $p$ (and set $s(A)=-\infty$ if no such $s$ exists).

Demailly (\cite{De} Lemma 8.6) proves that given any $s\ge0$, if $H^0(kA)$
generates $(k(n+s)+1)$-jets at $p$, then $H^0(A+K_X)$ generates $s$-jets
at $p$. Setting
$s=(n+1)(m-1)$ (and noting that $k(n+1)m\ge k(n+s)+1$) we see that for $m\ge
1$ and any $k\ge1$,
$$
s(kA) \ge k(n+1)m \quad\Rightarrow\quad s(A+K_X)\ge (n+1)(m-1).
$$
Applying this again to $A+K_X$ (with $k=1$) implies that $s(A+2K_X)\ge(n+1)(m-2)$
and so inductively $s(A+rK_X)\ge(n+1)(m-r)$ for $r\le m$. Setting $r=m$,
\begin{equation}\label{generationofadjointbundles}
s(kA) \ge k(n+1)m \quad\Rightarrow\quad s(A+mK_X)\ge0.
\end{equation}

Fix any rational number $m/M\in(0,\epsilon(p,L))$, so that
$\epsilon(p,(n+1)ML)>(n+1)m$. Then by the identity (\cite{De} Lemma 7.6)
$$
\epsilon(p,A) = \lim_{k\to\infty}\frac{s(kA)}k\,,
$$
applied to $A=(n+1)ML$,
we can find a $k\gg0$ such that $s(k(n+1)ML)>k(n+1)m$. Thus by
(\ref{generationofadjointbundles}), $(n+1)ML+mK_X$ has a section (not vanishing
at $p$), and so $((n+1)ML+mK_X).L^{n-1}\ge0$. Therefore
$((n+1)L+\epsilon(p,L)K_X).L^{n-1}\ge0$.
\end{proof}

\begin{rmk}
If $X=\PP^n$ we may as well assume that $L$ is the
hyperplane bundle as slope stability is invariant under rescaling. So
$\epsilon(p)=1$ for any point $p$ and $-K=(n+1)L$, which gives
equality in (\ref{pts}), so $\PP^n$ is at best K-semistable (since
$L^k\otimes\I_p^k$ is generated by global sections).
In fact $\PP^n$ is K-polystable in the sense of Definition \ref{defstab},
and the deformation to the normal cone of $\{p\}$ has semi-ample
$\L_c$ for $c=1$, pulled back from a blow down to $\PP^n\times\C$.
So the degeneration is a product family with a nontrivial
$\C^\times$-action, and the semistability we are seeing comes from a
$\C^\times$-action on $\PP^n$ with Donaldson-Futaki invariant zero.
\end{rmk}

\subsection*{Asymptotic slope Chow stability}

If instead of K-instability we are interested in asymptotic stability,
then the relevant slope function is more complicated.

Let $(X,L)$ be a polarised manifold, $r$ be a positive integer, and 
$Z$ a subscheme of $X$ such that on the blow up $p\colon\widehat X\to X$
of $X$ along $Z$,
\begin{equation} \label{condition}
\mathbf Rp_*\O(-jE)=p_*\O(-jE)=\I_{\!Z}^j\quad\forall j\ge0.
\end{equation}
(E.g. if $Z\subset X$ are both smooth in a neighbourhood of $Z$.)
Then for \emph{all} $k$ and $x$ (with $k,\,xk\in\mathbb N$) we have the exact
formula
$$
\chi(L^k\otimes\I_{\!Z}^{xk}) = a_0(x)k^n + a_1(x)k^{n-1} +\ldots+a_n(x).
$$
Letting $B_i$
denote the Bernoulli numbers, define $\beta_0=1,\,\beta_1=\frac12$ and
$\beta_i=\frac{B_i}{i!}$ for $i\ge 2$. Then we define the \emph{asymptotic
Chow slope} $\eta_c(\I_Z)$ of $\I_Z$ for $c\le\epsilon(Z)$ to be
\begin{eqnarray*}
\eta_c(\I_Z)&=&\eta_c(\I_Z,L,r) = \sum_{j=0}^{n+1} c_j r^{n+1-j} \\
&=& r^{n+1} + \mu_c(\I_Z)r^n + \ldots + c_{n+1},
\end{eqnarray*} 
where
$$
c_j = \frac{\sum_{i=0}^j
\int_0^c\beta_ia_{j-i}^{(i)}(x) dx}{\int_0^c a_0(x)dx}\,.
$$
I.e. $c_0=1,\ c_1=\frac{\int_0^c a_1(x)+\frac{a_0'(x)}2dx}{\int_0^c a_0(x)dx}
=\mu_c(\I_Z)$,
$$
c_2=\frac{1}{\int_0^c a_0(x) dx} \int_0^c
a_2(x) + \frac{a_1'(x)}{2} + \frac{a_0^{''}(x)}{12}\,dx, \quad\mathrm{etc.}
$$
Defining the asymptotic Chow slope of $X$ to be the slope of the empty
subscheme,
$$
\eta\_X(r) = \eta(\O_X,L,r) = \frac{r\chi(L^r)}{a_0} = r^{n+1} + \mu(X)
r^n + \frac{a_2}{a_0}r^{n-1} + \ldots + \frac{a_n}{a_0}r,
$$
we say (with great difficulty) that $X$ is asymptotically
Chow slope strictly destabilised by $Z$ if, for all $r\gg0$,
$$
\eta_c(\I_Z,r)>\eta\_X(r).
$$

\begin{thm} \label{chow}
  If a polarised variety is asymptotically Chow slope strictly
  destabilised by $Z$ satisfying (\ref{condition}) then it is
  asymptotically Chow strictly unstable.  
\end{thm}

\begin{proof}
  The proof is almost the same as for Theorem
  \ref{thm:kstableslopestable}.  We calculate the coefficient
  $e_{n+1}(r)$ of $k^{n+1}$ in $\tilde{w}_{r,k}$
  (\ref{normalised}),
$$
e_{n+1}(r)\ =\ b_0r\chi(L^r)-a_0w(r)\ =\ (b_0+ca_0)r\chi(L^r)-
a_0\big(w(r)+cr\chi(L^r)\big)
$$
(instead of its $r^n$ coefficient the Donaldson-Futaki invariant) for the
deformation to the normal cone. By Theorem \ref{thm:instability}, if
$e_{n+1}(r)<0$ then $(X,L)$ is Chow unstable with respect to $r$. By the
continuity of $\eta_c$ we may take $c<\epsilon(Z)$,
so that $w(r)$ is calculated by Theorem \ref{degen} and $b_0$ by Proposition
\ref{integrals}, giving
$$
e_{n+1}(r)<0 \quad\Longrightarrow\quad
\int_0^ca_0(x)dx\,r\chi(L^r)<a_0\sum_{j=1}^{cr}\chi(L^r\otimes\I_{\!Z}^j).
$$
Instead of the trapezium rule (\ref{easyeuler}) we use the fact \cite{Hi}
that for any polynomial $f$,
$$
\sum_{j=1}^{cr} f(j/r) = \int_0^c \sum_{i=0}^n
\beta_i \frac{f^{(i)}(x)}{r^{i-1}}dx.
$$
(The case $f(x)=x^m$ follows from the definition of $\beta_i$, implying the
general case by linearity.) The theorem follows
by applying this  to the polynomial $f(x)= a_0(x)r^n + a_1(x)r^{n-1}
+\ldots+a_n(x)=\chi(L^r\otimes\I_{\!Z}^{xr})$. 
\end{proof}

\section{Simplifying arbitrary test configurations} \label{testconfigs}

We begin with an important technical result allowing us to calculate weights
of a test configuration $(\Y,\O_\Y(k))$ in terms of a semi test configuration
$(\X,\L^k)$ that dominates it.

So let $(\Y,\O_\Y(1))$ be an equivariantly polarised flat $\C^\times$-family
with general fibre $(X,L)$, and fix another flat $\C^\times$-family $(\X,\L)\to\C$
with a birational $\C^\times$-equivariant map
$$
p\colon(\X,\L)\to(\Y,\O_\Y(1)) \qquad \text{such that }\ \L=p^*\O_\Y(1)
\ (\text{equivariantly}).
$$

\begin{prop} \label{stein}
If $X$ is normal then there exists an $a\ge0$ such that
$$
w(H^0_{\Y_0}\!(\O_\Y(k)))
=w\left(H^0_\X(\L^k)\big/tH^0_\X(\L^k)\right)-ak^n+O(k^{n-1}).
$$
\end{prop}

\begin{proof}
$p$ factors through the its Stein factorisation (\cite{Ha} Corollary
III.11.5) as
$$
(\X,\L)\Rt{q}(\X'=\Proj\bigoplus_kH^0_\X(\L^k),\O\_{\!\X'}(1))\Rt{p'}(\Y,\O_\Y(1)).
$$
$q_*$ induces an equivariant isomorphism between $H^0_\X(\L^k)$ and
$H^0_{\X'}(\O\_{\!\X'}(k))$ for $k\gg0$, and $\X'$ is flat over $\C$ since
$H^0_\X(\L^k)$ has no $t$-torsion by the flatness of $\X$.
Thus replacing $\X$ by $\X'$ and $p$ by $p'$ if necessary, we may assume
that $p$ is finite and so $\L=p^*\O_\Y(1)$ is ample. Moreover, the general
fibre of $\X'$ is $(X,L)$ by the normality of $X$, so
$(\X,\L^k)$ is now a genuine test configuration for $(X,L^k)$.

So we have a $\C^\times$-equivariant morphism of polarised families
$(\X\rt{\Pi\,}\C,\L)\Rt{p}(\Y\rt{\pi}\C,\O_\Y(1))$, and we wish to relate
the total weights of the $\C^\times$-actions on
$H^0_{\Y_0}(\O_\Y(k))$ and $H^0_\X(\L^k)\big/tH^0_\X(\L^k)$; but the latter
is now isomorphic to $H^0_{\X_0}(\L^k)$ for $k\gg0$ by ampleness and flatness.
$p^*$ induces an exact sequence
$$
0\to\O_\Y\to p_*\O_\X\to Q\to0,
$$
for some cokernel $Q$ supported on $\Y_0$ (since $p$ is an isomorphism away
from the central fibres).
So $t^sQ=0$ for some $s\ge0$. We first give the argument for $s=1$
to illustrate the more technical general case.

Tensoring with $\O_\Y(k)$ and pushing down to $\C$ gives an exact sequence on
$\C$,
$$
0\to\pi_*(\O_\Y(k))\to\Pi_*(\L^k)\to Q_k\to0,
$$
where $Q_k=\pi_*(Q\otimes \O_\Y(k))$ satisfies $tQ_k=0$. Therefore its restriction to $0\in\C$
is also isomorphic to $Q_k$, and Tor$_1(Q_k,\O_0)=\ker\big(Q_k\Rt{\times t}Q_k\big)\cong
Q_k$, giving the exact diagram
$$
\begin{array}{ccccccccc}
& 0 & \to & \pi_*(\O_\Y(k)) & \Rt{p^*} & \Pi_*(\L^k) & \to & \!\!Q_k & \to0
\\ & \downarrow && \downarrow && \downarrow && \downarrow\wr \\
0\to & Q_k\! & \to & H^0_{\Y_0}(\O_\Y(k)) & \Rt{p^*} & H^0_{\X_0}(\L^k)
& \to & Q_k|\_0\! & \to0 \\
&&& \downarrow && \downarrow \\ &&& 0 && 0
\end{array}
$$
where the vertical arrows are restriction to $0\in\C$, and the flatness of
$\pi,\,\Pi$ and ampleness of $\O_\Y(1),\,\L$ give the two
central terms and ensure that Tor$_1(\Pi_*(\L^k),\O_0)=0$.

These are maps of $\C^\times$-modules, except that the left hand $Q_k$
has weight shifted by $-1$ (i.e.\ is isomorphic to $Q_k\otimes\langle t\rangle$
as a $\C^\times$-module). We see this by exhibiting an explicit weight-$(-1)$
isomorphism $\delta$ from $Q_k$ to the kernel of the lower $p^*$. Given $q\in
Q_k$, choose a lift $\hat q\in\Pi_*(\L^k)$. $t\hat q$ has zero image in $Q_k$
so is in the image of some $f\in\pi_*\O_\Y(k)$. $f|\_0\in H^0_{\Y_0}(\O_\Y(k))$
maps under $p^*$ to $t\hat q|\_0\in  H^0_{\X_0}(\L^k)$, i.e.\ zero, so is
an element $\delta(q)\in\ker p^*=\,$Tor$_1(Q_k,\O_0)$. Since this map involved
multiplication by $t$ it has weight $-1$.

So in this case of $s=1$, $w(H^0_{\Y_0}\!(\O_\Y(k)))=w(H^0_{\X_0}(\L^k))-\dim Q_k$, and
$\dim Q_k=ak^n+O(k^{n-1})$, where $a=$\,rank$(Q_k)\O_\Y(1)^n/n!\ge0$ since $Q_k$ is
supported on the $n$-dimensional central fibre.

For general $s$ we still have the exact sequence
$$
0\to\mathrm{Tor}_1(Q_k,\O_0)\to H^0_{\Y_0}(\O_\Y(k))\Rt{p^*} H^0_{\X_0}(\L^k)
\to Q_k|\_0\to0,
$$
where again the first term $\ker\big(Q_k\Rt{\times t}Q_k\big)$ is isomorphic
to the last $Q_k|\_0$. For instance if we pick an isomorphism $Q_k\cong\bigoplus_i
\C[t]/(t^{j_i})$, each $j_i\le s$, then $Q_k|\_0\cong\bigoplus_i\C[t]/(t)$
maps isomorphically via $\oplus_it^{j_i-1}$ to $\bigoplus_i(t^{j_i-1}\C[t])
/(t^{j_i})=\ker\big(Q_k\Rt{\times t}Q_k\big)$. More invariantly, filter
$Q^k$ by $\ker t^j\subseteq\ker t^{j+1}\subseteq\ldots\subseteq\ker t^s=Q^k$,
and so $Q_k|\_0$ by $\ker t^j/(\im t\cap\ker t^j)$ with graded pieces $V^j=
\ker t^j/(\ker t^{j-1}+\im t\cap\ker t^j)$. Then there exists a weight-$(-j)$
``multiplication
by $t^j$" map taking $V^j$ into $H^0_{\Y_0}(\O_\Y(k))$ whose sum over $j\le
s$ gives an isomorphism to Tor$_1(Q_k,\O_0)=\ker p^*$.

The map is defined roughly as before. Given $q\in V^j$, lift to $\ker(Q_k\rt{t^j}Q_k)$
and then to $\hat q\in\Pi_*(\L^k)$. The image of $t^j\hat q$ in $Q_k$
is zero by construction, so is in the image of $f\in\pi_*\O_\Y(k)$. $f|\_0\in
H^0_{\Y_0}(\O_\Y(k))$ maps to $t^j\hat q|\_0\in H^0_{\X_0}(\L^k)$, i.e.\ to
zero since $j\ge1$. Thus $f|\_0=\delta(q)$ for some $\delta(q)\in\ker p^*$.
$\ker t^{j-1}$ and $\im t\cap\ker t^j$ clearly map to zero via this construction,
so $\delta$ is well defined. The upshot is that
$$
w(H^0_{\Y_0}\!(\O_\Y(k)))=w(H^0_{\X_0}(\L^k))-\sum_{j=0}^sj\dim V_j.
$$
$\sum_{j=0}^sj\dim V_j$ is a polynomial in $k$ (since both weights $w$
are), so is $ak^n+O(k^{n-1})$ for some $a\ge0$ because $0\le\sum\dim V_j=\dim
Q_k=O(k^n)$ as before.
\end{proof}

\begin{rmk} \label{assumenormal}
In particular if $X$ is normal then for K-stability we need only
consider \emph{normal} test configurations. This is because
any test configuration $(\Y,\O_\Y(1))$ is dominated by its normalisation,
which also has general fibre $X$ if $X$ is normal. The pullback of $\O_\Y(1)$
is ample, so some twist gives another test configuration which is \emph{less
stable} than $(\Y,\O_\Y(1))$, in the sense that it has the same weight to
leading order and a smaller Donaldson-Futaki invariant, by Proposition
\ref{stein} above. 

Similarly for Chow stability we may compute the Chow weight of any test
configuration on its normalisation. Of course most of our test configurations
have nonnormal central fibre, however.
\end{rmk}

Given any test configuration $(\Y,\O_\Y(1))$ (\ref{any}) we now build
inductively the semi test configuration $(\B_{I_r}(X\times\C),L-E)$
(\ref{ideals}) that dominates it, starting from the deformation to the normal
cone $(\B_{I_1}(X\times\C),L-P)$ of $Z_0$.

So let $\X^1\rt{\pi^1}X\times\C$ denote $\B_{I_1}(X\times\C)$, i.e.\ the blow
up in $Z_0\times\{0\}$:
$$
\X^1=\Proj\bigoplus_k\,\S_k, \qquad\qquad \S_k:=(\I_0+(t))^k.
$$
Recall that the central fibre of $\X^1$ is $\widehat X\cup_eP$, where $\widehat
X$ is the blow up of $X$ along $Z_0$ with exceptional divisor $e$, and $P$
is the exceptional divisor of $\pi^1$.
$P$ is a projective cone over $Z_0$ (Proj of the graded algebra over $Z_0$
with $k$th graded piece $\bigoplus_{i=0}^{k-1}\I_0^i/\I_0^{i+1}$) -- the projective
completion of the normal cone to $Z_0\subset X$. Its zero
section $Z_0'$ is a copy of $Z_0$ which fits into a flat
family with the $Z_0\times\{t\}$ in each fibre $X_t$, which we see as follows.

The proper transform $\overline{Z_0\times\C}$ is defined
by the graded sheaf of ideals generated by $\I_0\subset\S_1=\I_0+(t)$
in the graded sheaf of algebras $\bigoplus_k\S_k$. That is, $\O(-P)\otimes
\I_{\overline{\!Z_0\times\C}}$ is generated by the sections of $\I_0\subset\S_1$.
It is abstractly isomorphic
to the blow up of $Z_0\times\C$ along its intersection with $Z_0\times\{0\}$,
but this is a divisor in $Z_0\times\C$, so $\overline{Z_0\times\C}\cong
Z_0\times\C$. The central fibre $Z_0'\cong Z_0$ is
defined by the graded sheaf of ideals generated by $\I_0+t\S_1=\I_0+(t^2)
\subset\S_1$. 

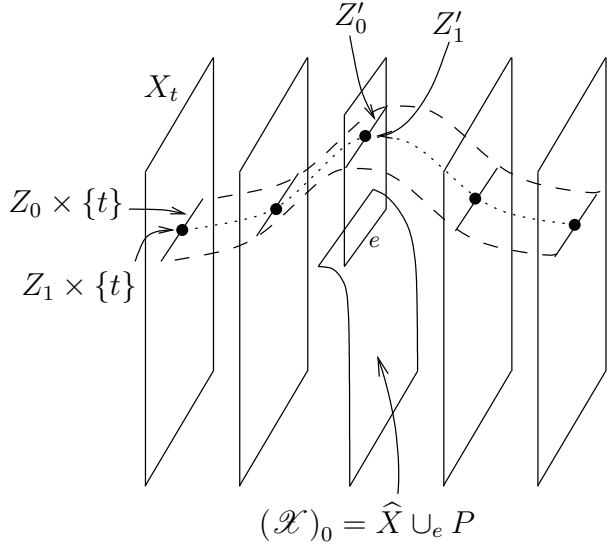
\begin{figure}[h]
\center{\input{normalcone.pstex_t}}
\end{figure}

Similarly the proper transform of $Z_1\times\C$ is the blow up of $Z_1\times\C$
along its intersection $Z_1\times\{0\}$ with $Z_0\times\{0\}$; that is 
$\overline{Z_1\times\C}\cong Z_1\times\C$. It is defined by the graded sheaf
of ideals generated by $\I_0+\I_1\S_1=\I_0+t\I_1\subset\S_1$,
with central fibre $Z_1'\subseteq Z_0'\subset P$ isomorphic to $Z_1$ and
defined by
\begin{equation} \label{I1'}
\I_0+(\I_1+(t))\S_1=\I_0+t\I_1+(t^2)\ \subset\ \S_1.
\end{equation}

We now form $\X^2$ by blowing up $\X^1$ in $Z_1'$. Since $\I_0+t\I_1+(t^2)$
is just $I_2$ (\ref{ideals}), we have basically shown that $\X^2$ dominates
$\B_{I_2}(X\times\C)$. Precisely, we have maps $\X^2\rt{\pi^2}
\X^1\rt{\pi^1}X\times\C$, and set $E_1:=(\pi^1)^*P,\ E_2$ to be the exceptional
divisor of $\pi^2$, and $E$ to be the exceptional divisor of $\B_{I_2}(X\times\C)$.
While $\O(-cE_1-eE_2)$ is relatively
ample for $0<e<c$, it is only semi-ample for $e=c$\,:

\begin{prop} \label{p}
$\X^2\to X\times\C$ factors through a map $p^2\colon\X^2\to\B_{I_2}(X\times\C)$.
Under this map, $(p^2)^*(\O(-E))=\O(-E_2-E_1)$.
\end{prop}

\begin{proof}
For $k$ sufficiently large, $\pi^1_*\O(-kE_1)=(\I_0+(t))^k=\S_k$, in which
the ideal $(\I_0+t\I_1+(t^2))^k$ defines the $k$th power of the ideal of
$Z_1'$ (\ref{I1'}). That is, $(\pi^1\comp\pi^2)^*(\I_0+t\I_1+(t^2))^k\cong
(\pi^2)^*\I_{\!Z_1'}^k(-kP)=\O(-kE_1-kE_2)$.
Therefore the sections $(\I_0+t\I_1+(t^2))^k\subset(\pi^1\comp\pi^2)_*\O(-kE_1-kE_2)$
generate $\O(-kE_1-kE_2)$ and so define a regular map from $\X^2$ to
$\Proj\bigoplus_{k\gg0}(\I_0+t\I_1+(t^2))^k$, i.e.\ to
$$
\Proj\bigoplus_{k\gg0}I^k_2=\B_{I_2}(X\times\C). \vspace{-9mm}
$$ \vspace{2mm}
\end{proof}

In fact the contraction $p^2$ just collapses the restriction $P\res{Z_1}$
of the cone $P\to Z_0$ down to $Z_1$, but we will not need this.

Similarly in $\X^2$ there is a copy $Z_1''$ of $Z_1$, sitting in a flat
family with $Z_1\subset X$ as the central fibre of the proper transform
$\overline{Z_1\times\C}$. In the coordinate ring
$\bigoplus_k\S_k=\bigoplus_k(\I_0+t\I_1+(t^2))^k$ (we are recycling the symbol
$\S_k$ and shall do so again below) pulled back from $\X$ by the above map
$p^2$ (\ref{p}), $\I_{\overline{Z_1\times\C}}$ is generated by $\I_0+t\I_1\subset
(\I_0+t\I_1+(t^2))$, since this is the largest ideal that localises to $\I_1$
when we invert $t$ and work on $X\times\C^\times$. Thus $\I_{Z_1''}$ is generated
by $\I_0+t\I_1+t(\I_0+t\I_1+(t^2))=\I_0+t\I_1+(t^3)$. 

So there is also a $Z_2''\subset Z_1''$, isomorphic to $Z_2$, inside it.
$\overline{Z_2\times\C}$ has ideal generated by $\I_0+t\I_1+t^2\I_2\subset
\S_1=\I_0+t\I_1+(t^2)$, since this is the largest ideal that localises to $\I_2$
when we invert $t$ and work on $X\times\C^\times$. Thus its central fibre
$Z_2''$ is defined by the ideal
\begin{equation} \label{I3}
\I_0+t\I_1+t^2\I_2+t\S_1=\I_0+t\I_1+t^2\I_2+(t^3).
\end{equation}

Blowing it up gives $\X^3$, with exceptional divisor $E_3$ (and we denote
the pullbacks to $\X^3$ of $E_1,\,E_2$ by the same notation). Then, just
as in Proposition \ref{p}, the pushdown of $\O_{\X^3}(-E_3-E_2-E_1)$ to $X\times\C$
is generated by $I_3=\I_0+t\I_1+t^2\I_2+(t^3)$ by (\ref{I3}), so $\X^3\to X\times\C$
factors through $\B_{I_3}(X\times\C)$.

Inductively we obtain $\X^s\to X\times\C$
as the blow up $\pi^s$ along $Z_{s-1}^{(s-1)}\subset\X^{s-1}$,
the central fibre of the proper transform of $Z_{s-1}
\times\C$. The coordinate ring of $\X^{s-1}$ over $X\times\C$ has $k$th graded
piece $\S_k=(\I_0+t\I_1+\ldots+t^{s-2}\I_{s-2}+(t^{s-1}))^k$, and the ideal
of the proper transform of $\overline{Z_{s-1}\times\C}$ is generated by
$\I_0+\ldots+t^{s-2}\I_{s-2}+t^{s-1}\I_{s-1}\subset\S_1=\I_0+\ldots+t^{s-2}\I_{s-2}
+(t^{s-1})$, since this is the largest ideal that localises to $\I_{s-1}$
when we invert $t$ and work on $X\times\C^\times$. Therefore the ideal of
its central fibre $Z_{s-1}^{(s-1)}$ is generated by
$$
\I_0+\ldots+t^{s-2}\I_{s-2}+t^{s-1}\I_{s-1}+t\S_1\ =\ 
\I_0+\ldots+t^{s-2}\I_{s-2}+t^{s-1}\I_{s-1}+(t^s).
$$
Thus the pushdown of $\O\_{\!\X^s}(-E_s-\ldots-E_1)$ contains $\I_0+t\I_1+
\ldots+t^{s-1}\I_{s-1}+(t^s)$, giving the following, as in Proposition \ref{p}.

\begin{thm} \label{pr}
$\X^s\to X\times\C$ factors through a map $p^s\colon\X^s\to\B_{I_s}(X\times\C)$,
where $I_s=\I_0+t\I_1+\ldots+t^s\I_{s-1}+(t^s)$ (\ref{ideals}).
Under this map, $(p^s)^*(\O(-E))=\O(-E_s-\ldots-E_1)$. \hfill$\square$
\end{thm}

In turn any test configuration $(\Y,\O_\Y(1))$ is dominated by a map $\phi$
from some $\B_{I_r}(X\times\C)$ (\ref{any}), giving $\rho^r:=\phi\comp p^r\colon
\X^r\to\Y$. Denote by $L_r$ the semi-ample line bundle $(\pi^1\comp\ldots\comp
\pi^r)^*L-E_r-\ldots-E_1$ on $\X^r$, so that $(\rho^r)^*\O_\Y(1)=L_r$.
Then by the above and Proposition \ref{stein} we have

\begin{cor} \label{ref}
The total weight on (the $k$th twist of) an arbitrary test configuration
$(\Y,\O_\Y(1))$ can be calculated on $(\X^r,L_r)$ by
$$
w(H^0_{\Y_0}(\O_\Y(k))=
w\left(H^0_{\X^r}(L_r^k)\big/tH^0_{\X^r}(L_r^k)\right)-ak^n+O(k^{n-1}).
$$
\end{cor}

\begin{cor} \label{seshadriIr}
The Seshadri constant $\epsilon(I_r)$ of the ideal $I_r$ (\ref{ideals})
is $\min\{\epsilon(Z_i)\}_{i=1}^{r-1}$.
\end{cor}

\begin{proof}
If $c\in(0,\min\{\epsilon(Z\_i)\}_{i=1}^{r-1})\cap\Q$ then for $k\gg0$, the sections
$$
H^0_X(L^k\otimes\I_0^{ck})\,\oplus\,t^{ck}H^0_X(L^k\otimes\I_1^{ck})\,\oplus\,
\ldots\,\oplus\,t^{(r-1)ck}H^0_X(L^k\otimes\I_{r-1}^{ck})\,\oplus\,t^{rck}H^0_X(L^k)
$$
saturate $I^{ck}_r$ (\ref{ideals}). Therefore the Seshadri constant of $I_r$
(which is $\ge1$ by the semi ampleness of $L-E$ in (\ref{any})) is at least
$\min\{\epsilon(Z\_i)\}_{i=1}^{r-1}$.

We prove the opposite inequality inductively. The induction begins for
$r=0$ by Proposition \ref{ample}; suppose it is true for $r$. That is,
$L-xE$ is nef on $\B_{I_r}(X\times\C))$ precisely when
$x\in[0,\min\{\epsilon(Z_i)\}_{i=1}^{r-1}]$. Equivalently, pulling back by
$p^r$, $L-x(E_1+\ldots+E_r)$ is nef on $\X^r$ if and only if $x\in[0,
\min\{\epsilon(Z_i)\}_{i=1}^{r-1}]$. Pick $c$ such that
$L-c(E_1+\ldots+E_{r+1})$ is nef on $\X^{r+1}$; we must show that
$c\le\min\{\epsilon(Z_i)\}_{i=1}^{r-1}$ and $c\le\epsilon(Z_r)$.

To show the former, we claim that since $L-c(E_1+\ldots+E_{r+1})$
is nef on $\X^{r+1}$, $L-c(E_1+\ldots+E_r)$
is nef on $\X^r$. Given any irreducible proper curve $C\subset\X^r$ not entirely
contained in $Z_r^{(r)}$, let $\overline C$ denote its proper transform in
$\X^{r+1}$. Then $[L-c(E_1+\ldots+E_r)].C\ge[L-c(E_1+\ldots+E_{r+1})].\overline
C\ge0$ since $\overline C.E_{r+1}\ge0$. On the other hand if $C\subset Z_r^{(r)}$,
then there is an isomorphic copy $C'\subset Z_r^{(r+1)}\subset\X^{r+1}$ such
that $[L-c(E_1+\ldots+E_r)].C=[L-c(E_1+\ldots+E_{r+1})].C'\ge0$. So indeed
$L-c(E_1+\ldots+E_r)$ is nef on $\X^r$ and so $c\le\min\{\epsilon(Z_i)\}_{i=1}^{r-1}$
by induction.

Secondly, fix $c$ such that $L-c(E_1+\ldots+E_r)$
is nef on $\X^r$. $Z_r^{(r)}\cong Z_r$ lies in the central fibre $(\X^r)\_0$,
fitting into a flat family with
$Z_r\times\C^\times$ away from the central fibre. Seshadri constants are
lower semicontinuous in polarised families, so the Seshadri constant (with
respect to $L-c(E_1+\ldots+E_r)$) of $Z_r^{(r)}$ \emph{inside
the central fibre $(\X^r)\_0$} is at most $\epsilon(Z_r)$ (with
respect to $L$, since this is $L-c(E_1+\ldots+E_r)$ restricted to a general
fibre). The Seshadri constant of $Z_r^{(r)}$ inside the whole of
$\X^r$ can only be smaller still; therefore if $L-c(E_1+\ldots+E_r)-cE_{r+1}$
is nef then $c\le\epsilon(Z_r)$.
\end{proof}

We could contract each $\X^s$ using $L_s$; this would give an
isomorphism in a neighbourhood of $Z_s^{(s)}$ by (\cite{Ha}
Proposition II.7.3) since $L_s|_{Z_s^{(s)}}
\cong L|\_{Z_s}$ is ample \emph{and} $L_s^k\otimes\I_{\!Z_s^{(s)}}$
is globally generated (by sections of $L^k\otimes(\I_0+t\I_1+\ldots+
t^s\I_s+(t^{s+1}))$ on $X$). Thus we could proceed inductively with these
contracted $\X^s$s with \emph{ample} line bundles on them, but since we
cannot currently seem to get significantly better estimates by working with
ample bundles we proceed with the semi-ample $(\X^s,L_s)$ and contract at
the last, $r$th stage. By Corollary \ref{ref} we lose nothing by ignoring
this contraction and simply calculating the weight on $\X^r$.

So, modulo the contraction, we have exhibited any test
configuration (\ref{any}) as a finite number of blow ups (starting with
$X\times\C$) \emph{in subschemes $Z^{(i)}_i$ supported in the scheme theoretic
central fibre} that themselves sit in flat families with the $Z_i\subset X$. We
calculate the weight on such a blow up in Theorem \ref{big}, for which
we need two preliminary results.

The following Proposition is the appropriate generalisation to general
test configurations of the case $Z\times\C\subset X\times\C$ used in (\ref{St}),
(\ref{qt}) and (\ref{error}). We will apply it to the flat families
$\Z=\overline{Z_s\times\C}\,\subset\X^s\to\C$
and their thickenings $k\Z$, when these thickenings are also flat.

\begin{prop} \label{central}
Fix flat families $\Z\subset\X\to\C$ with central fibres
$Z'\subset\X_0$, such that the thickenings $k\Z\subset\X\to\C$
are also flat over $\C$. Let $\I\_{\!Z'\subset\X}$ (respectively
$\I_{Z'}$) denote the ideal sheaf of $Z'\subset\X$ ($Z'\subset\X_0$).
Then
$$
\frac{\I_{\!Z'\subset\X}^k}{t\I_{\!Z'\subset\X}^k}\ \cong\
\I_{\!Z'}^k\ \oplus\ \bigoplus_{j=1}^k\,t^j
\frac{\I_{\!Z'}^{k-j}}{\I_{\!Z'}^{k-j+1}}\,.
$$
\end{prop}

\begin{proof}
Flatness of $\X\to\C$ and $j\Z\to\C$ imply the exactness of the
bottom two rows of the following, from which follows the exactness of the
top row.
$$
\begin{array}{ccccccc}
&0&&0&&0 \\
& \downarrow && \downarrow && \downarrow \\
0\to & \I_\Z^j & \rt{t} & \I_\Z^j & \to & \I_{\!Z'}^j & \to0 \\
& \downarrow && \downarrow && \downarrow \\
0\to & \O_\X & \rt{t} & \O_\X & \to & \O_{\X_0} & \to0 \\
& \downarrow && \downarrow && \downarrow \\
0\to & \O_{j\Z} & \rt{t} & \O_{j\Z} & \to & \O_{jZ'} & \to0 \\
& \downarrow && \downarrow && \downarrow \\
&0&&0&&0
\end{array}
$$
Chasing through either the top two rows or the first two columns then shows
that in $\O_\X$, $(t)\cap\I_\Z^j=t\I_\Z^j$. Also
by flatness there is a similar diagram with $t$ replaced by $t^i$ (and the
right hand column suitably modified) showing that in fact
\begin{equation} \label{one}
(t^i)\cap\I_\Z^j=t^i\I_\Z^j.
\end{equation}
Applying $H^0_\X(L\otimes\,\cdot\,)$ gives a similar exact diagram without
the right hand and lower zeros, so the same argument shows that
\begin{equation} \label{one'}
(t^i)\cap H^0_\X(L\otimes\I_\Z^j)=t^iH^0_\X(L\otimes\I_\Z^j).
\end{equation}
The top row also gives $\I_\Z^j/t\I_\Z^j=\I_{\!Z'}^j$, which implies that
\begin{equation} \label{four}
\frac{\I_\Z^j}{\I_\Z^{j+1}+t\I_\Z^j}=\frac{\I_{\!Z'}^j}{\I_{\!Z'}^{j+1}}\,.
\end{equation}
$\I\_{Z'\subset\X}=\I_\Z+(t)$, so $\I_{\!Z'\subset\X}^k=
\sum_{j=0}^k t^j\I_\Z^{k-j}$ and
\begin{equation} \label{five}
\frac{\I_{\!Z'\subset\X}^k}{t\I_{\!Z'\subset\X}^k}=
\frac{\sum_{j=0}^k t^j\I_\Z^{k-j}}{\sum_{j=0}^k t^{j+1}\I_\Z^{k-j}}=
\sum_{j=0}^k\frac{t^j\I_\Z^{k-j}}{t^j\I_\Z^{k-j}\ \cap\ \sum_{i=0}^k
t^{i+1}\I_\Z^{k-i}}\,,
\end{equation}
where we use the fact that in an abelian category, for $A,B,C\subset V$
with $C\subset A+B$, we have $(A+B)/C=A/(A\cap C)+B/(B\cap C)$ in $V/C$.
We claim that
\begin{equation} \label{three}
t^j\I_\Z^{k-j}\ \cap\ \sum_{i=0}^kt^{i+1}\I_\Z^{k-i}=
t^j\I_\Z^{k-j+1}+t^{j+1}\I_\Z^{k-j},
\end{equation}
except for $j=0$ when the right hand side becomes $t\I_\Z^k$.
The inclusion $\supseteq$ is clear. For $\subseteq$, consider the left
hand side:
\begin{eqnarray*}
t^j\I_\Z^{k-j} \!\!&\cap&\!\! \left((t\I_\Z^k+\ldots+t^j\I_\Z^{k-j+1})+
(t^{j+1}\I_\Z^{k-j}+\ldots+(t^{k+1}))\right) \\
\subseteq t^j\I_\Z^{k-j} \!\!&\cap&\!\! \left(t\I_\Z^{k-j+1}+(t^{j+1})\right).
\end{eqnarray*}
An element of this can be written as $t^jf=tg+t^{j+1}h$, that is
$t^{j-1}f-g=t^jh$, where $f\in\I_\Z^{k-j},\ g\in\I_\Z^{k-j+1}$
and so $t^{j-1}f-g\in\I_\Z^{k-j}$.
So $t^jh\in(t^j)\cap\I_\Z^{k-j}$, which by (\ref{one}) is $t^j\I_\Z^{k-j}$.
Thus $h$ may be taken to lie in $\I_\Z^{k-j}$. Similarly
$g\in(t^{j-1})\cap\I_\Z^{k-j+1}=t^{j-1}\I_\Z^{k-j+1}$.
Therefore $t^jf=tg+t^{j+1}h\in t^j\I_\Z^{k-j+1}+t^{j+1}\I_\Z^{k-j}$, proving
the inclusion.

So by (\ref{three}) and (\ref{four}), equation (\ref{five}) has become
$$
\frac{\I_{\!Z'\subset\X}^k}{t\I_{\!Z'\subset\X}^k}=\frac{\I_\Z^k}{t\I_\Z^k}+
\sum_{j=0}^k t^j \frac{\I_\Z^{k-j}}{\I_\Z^{k-j+1}+t\I_\Z^{k-j}}
=\I_{\!Z'}^k+\sum_{j=0}^k t^j\frac{\I_{\!Z'}^{k-j}}{\I_{\!Z'}^{k-j+1}}\,.
$$
To check that the sum is direct, intersect the $j$th numerator with the others:
$$
t^j\I_\Z^{k-j}\cap\sum_{p\ne j}t^p\I_\Z^{k-p}\subseteq t^j\I_\Z^{k-j}\cap
\left(\I_\Z^{k-j+1}+(t^{j+1})\right).
$$
By the same methods as before this lies in $t^j\I_\Z^{k-j+1}+t^{j+1}\I_\Z^{k-j}$,
which is the $j$th denominator, as required.
\end{proof}

To apply the above result inductively to the $\X^i$ requires the flatness
of the
thickenings $k(\overline{Z_i\times\C})$ in $\X^i$ of the proper transforms
of the $Z_i\times\C\subset X\times\C$. This is only automatic for $k=1$,
but also holds for arbitrary $k$ if the $Z_i$ and $X$ are all \emph{smooth},
or, we shall show, if:
\begin{eqnarray} \nonumber
\text{\emph{$X$ is reduced, each $Z_{r-1}\subseteq\ldots\subseteq Z_0$ is
a Cartier divisor in $X$, and any irreducible}} \hspace{-1cm} \\ \label{reduced}
\text{\emph{component common to any pair $Z_i,\,Z_j$ has the same
multiplicity in each.}} \quad
\end{eqnarray}
This odd looking condition is clearly satisfied if, for instance, each $Z_i$
is reduced.

\begin{prop} \label{flat}
Suppose that $X$ and the $Z_i$ satisfy (\ref{reduced}).
Then $j(\overline{Z_i\times\C})\subset \X^i$ is flat over $\C$ for each
$j\in\mathbb N$.
\end{prop}

\begin{proof}
Firstly, consider $\overline{j(Z_i\times\C)}$. Its ideal is defined by
those functions which, on restriction to $t\ne0$, lie in $\I_{\!Z_i\times
\C^\times}^j$. Since $t$ is invertible there, this implies that if $tf\in
\I\_{\overline{j(Z_i\times\C)}}$ then $f\in\I\_{\overline{j(Z_i\times\C)}}$\,.
Therefore the structure sheaf of
$\overline{j(Z_i\times\C)}$ has no $t$-torsion. It is also flat away from
$t=0$ (where it is $j(Z_i\times\C^\times)$); thus $\overline{j(Z_i\times\C)}$
is automatically flat over $\C$.

$\I\_{\overline{j(Z_i\times\C)}}\supseteq\I^j_{\overline{Z_i\times\C}}$;
therefore to prove that $j(\overline{Z_i\times\C})$ is flat over $\C$ it
is sufficient to prove the opposite inclusion
\begin{equation} \label{inclusion}
\I\_{\overline{j(Z_i\times\C)}}\ \subseteq\ \I^j_{\overline{Z_i\times\C}}\,.
\end{equation}

$\X^i$ is obtained by blowing up $\X^{i-1}$ in $Z_{i-1}^{(i-1)}$; we claim
the complement
$\X^i_\text{aff}$ of the proper transform $\overline{(\X^{i-1})\_0}$ of the
central fibre in $\X^i$ is affine over $X\times\C$ with coordinate ring
\vspace{3mm} \begin{equation} \vspace{-11mm} \label{coord} \end{equation}
\begin{multline*}
\!\!\sum_{a\_0,a\_1,\ldots,a\_{i-1}\ge0}\left(\!\frac{\I_0}{t^i}\!\right)^{\!a_0}
\!.\!\left(\!\frac{\I_1}{t^{i-1}}\!\right)^{\!a_1}
\!\ldots\!\left(\!\frac{\I_{i-1}}t\!\right)^{\!a\_{i-1}}\ = \\
\O_X\otimes\C[t]\ +\ \frac{\I_{\!i-1}}t\ +\ \frac{\I_{\!i-2}+\I_{\!i-1}^2}{t^2}
\ +\ \frac{\I_{\!i-3}+\I_{\!i-2}\I_{\!i-1}+\I_{\!i-1}^3}{t^3}\ +\ \ldots,
\end{multline*}
with its obvious ring structure. Each $\I_i$ is an $\O_X$-module, so the
whole ring inherits an $\O_X\otimes\C[t]$-module structure, corresponding
to the projection to $X\times\C$.

There are many ways to see this. One is to note that, away from $\overline
{(\X^{i-1})\_0}$\,, the map to $\B_{I_i}(X\times\C)$ of Theorem
\ref{pr} is an isomorphism. This is because, in the notation of that section,
sections of $\O(-E_1-\ldots-E_i)$ over $X\times\C$ do not contract the
exceptional divisor $E_i$ of the blow up of $\X^{i-1}$ in $Z_{i-1}^{(i-1)}$
(as noted in the remarks following the proof of Corollary \ref{seshadriIr}).
Let $E$ denote the exceptional divisor of
$\B_{I_i}(X\times\C)$, and $s_E\in H^0(\O(E))$ the canonical section
vanishing on $E$. For $k\gg0$, the sections of $\O(-kE)$ are the sections
of $I_i^k/s_E^k$ (that is, pull back sections of $I_i^k$ from
$X\times\C$ to $\B_{I_i}(X\times\C)$ and divide by $s_E^k$ to get a
regular section of $\O(-kE)$). $t^i\in I_i=\I_0+t\I_1+\ldots+
t^{i-1}\I_{i-1}+(t^i)$
defines the section $t^i/s_E$ of $\O(-E)$ which trivialises $\O(-E)$ over
$\X^i_\text{aff}$ -- the complement of its zero locus $\overline{(\X^{i-1})\_0}$.
Using this trivialisation identifies the functions on $\X^i_\text{aff}$ which
have poles of order $\le k$ on $\overline{(\X^{i-1})\_0}$ with
$$
\frac{I_i^k}{t^{ik}}\ =\ 
\left(\O+\frac{\I_{i-1}}{t}+\ldots+\frac{\I_0}{t^i}\right)^k;
$$
taking the limit as $k\to\infty$ gives the regular functions (\ref{coord}).

Alternatively, we can work inductively with the $\X^j$. A similar analysis
as above shows the coordinate ring of the complement of $\overline{\X_0}$
in the blow up of a family $\X$ over $\C$ in an ideal $I+(t)$ is
\begin{equation} \label{BI}
\O_\X+\frac It+\frac{I^2}{t^2}+\frac{I^3}{t^3}+\ldots
\end{equation}
Thus we find the coordinate ring of $\X^1_\text{aff}=\X^1\backslash\,
\overline{X\times\{0\}}$ is
$$
\O_X\otimes\C[t]\,+\,\frac{\I_0}t\,+\,\frac{\I_0^2}{t^2}\,+\,\frac{\I_0^3}{t^3}\,+\ldots
$$
Inside this the ideal of $\overline{Z_1\times\C}$ is $\I_1+\frac{\I_0}t+
\frac{\I_0^2}{t^2}+\ldots$ (as this is the largest ideal that localises
to $\I_1\otimes\C[t,t^{-1}]$ on $t\ne0$); applying (\ref{BI}) to this ideal
$I$ gives the coordinate ring of $\X^2_\text{aff}=\X^2\backslash\,
\overline{(\X^1)_0}$ as
$$
\O_X\otimes\C[t]\ +\ \frac{\I_1}t\ +\ \frac{\I_1^2+\I_0}{t^2}\ +\ 
\frac{\I_1^3+\I_0\I_1}{t^3}\ +\ \ldots
$$
But this is the $i=2$ case of (\ref{coord}), and inductively we recover it
for all $i$. In (\ref{coord}) we have the ideal
\begin{equation} \label{kZ}
\I\_{\overline{j(Z_i\times\C)}}\ =\ 
\sum_{a\_0,a\_1,\ldots,a\_{i-1}\ge0}\frac{\I_i^j\cap(\I_0^{a_0}.\I_1^{a_1}
\ldots\I_{i-1}^{a\_{i-1}})}{t^{ia_0}.t^{(i-1)a_1}\ldots t^{a\_{i-1}}}\,,
\end{equation}
as this is the largest ideal that localises to $\I_i^j\otimes\C[t,t^{-1}]$
on $t\ne0$. In the $j=1$ case, $\I_i\subseteq\I_0^{a_0}.\I_1^{a_1}
\ldots\I_{i-1}^{a\_{i-1}}$ unless $a_j=0\ \forall j$,
so $\I\_{\overline{Z_i\times\C}}$ differs from (\ref{coord}) only in the
first term $\I_i\subset\O_X\otimes\C[t]$:
\begin{equation} \label{Z}
\I\_{\overline{Z_i\times\C}}\ =\ \I_i\ +\!\sum_{a\_j\ge0,\ \sum_{j=0}^{i-1}a_j\ge1}
\frac{\I_0^{a_0}.\I_1^{a_1}
\ldots\I_{i-1}^{a\_{i-1}}}{t^{ia_0}.t^{(i-1)a_1}\ldots t^{a\_{i-1}}}\,.
\end{equation}
By (\ref{inclusion}) we are left with showing that each term of (\ref{kZ}) is
contained in the $j$th power of (\ref{Z}):
\begin{equation} \label{Zk}
\I^j_{\overline{Z_i\times\C}}\ =\sum_{p=0}^j\left(\sum_{a\_j\ge0,\
\sum_{j=0}^{i-1}a_j\ge p}\!\!\I_i^{j-p}\cdot\frac{\I_0^{a_0}.\I_1^{a_1}
\ldots\I_{i-1}^{a\_{i-1}}}{t^{ia_0}.t^{(i-1)a_1}\ldots t^{a\_{i-1}}}\right).
\end{equation}

We now work locally, where the conditions (\ref{reduced}) on $Z_i\subseteq
Z_{i-1}\subseteq\ldots\subseteq Z_0$ imply that
$\I_i=(f_i)$ and $\I_j=(g_jf_i),\ j\le i-1$ for some $g_j$
\emph{which do not divide} $f_i$. Therefore,
for any $a_j\ge0$ with $p:=\sum_{j=0}^{i-1}a_j\le j$, we have
\begin{multline*}
\I_i^j\cap(\I_0^{a_0}\ldots\I_{i-1}^{a\_{i-1}})=
(f_i^j)\cap(f_i^p.g_0^{a_0}\ldots g_{i-1}^{a\_{i-1}})
=(f_i^j.g_0^{a_0}\ldots g_{i-1}^{a\_{i-1}}) \\
=(f_i^{j-p}.f_i^p.g_0^{a_0}\ldots g_{i-1}^{a\_{i-1}})
=\I_i^{p-a}.\I_0^{a_0}\ldots\I_{i-1}^{a\_{i-1}},
\end{multline*}
using the fact the $\O_X$ is
torsion-free. This gives the desired inclusion.
\end{proof}

To apply Proposition \ref{flat} will involve replacing $X$ by a blow up on
which the pullbacks of the $Z_i$ are divisors. Pulling back the polarisation,
we find we are forced to work with a semi-ample line bundle. To this end,
for any $Z\subset X$ and \emph{semi-ample} $L\to X$, we define 
\begin{equation} \label{wkZ}
w_k(Z):=\sum_{j=1}^{ck}h^0_X(L^k\otimes\I_{\!Z}^j)-ckh^0_X(L^k)
\end{equation}
(cf. (\ref{degen}) and (\ref{ss})) for any $c\in\Q$ such that $L^k\otimes\I_{\!Z}^{ck}$
is saturated by its global sections. This equals
\begin{equation} \label{wkchi}
w_k(Z)=\sum_{j=1}^{ck}\chi\_X(L^k\otimes\I_{\!Z}^j)-ck\chi\_X(L^k)+O(k^n)
\qquad (O(k^{n-1})\text{ for $L$ ample}),
\end{equation}
by the following generalisation of Lemma \ref{error} applied to $Z$ and
$Z=\emptyset$.

\begin{lem} \label{wkerror}
For $Z\subset X,\ L\to X$ semi-ample, and $L^k\otimes\I_{\!Z}^{ck}$ saturated
by global sections for $k\gg0$,
$\sum_{j=1}^{ck}h^{\ge1}(L^k\otimes\I_{\!Z}^j)=O(k^n)$.
\end{lem}

\begin{proof}
Consider the
blow up of $X\times\PP^1$ in $Z\times\{0\}$ with exceptional divisor $P$
and line bundle $\O_{\PP^1}(ck)\boxtimes L^k-ckP$. This is semi-ample, generated
by $t^{ck}\boxtimes H^0_X(L^k)\,+\,s^{ck}\boxtimes H^0_X(L^k\otimes\I_{\!Z}^{ck})$,
where $s,\,t\in H^0(\O_{\PP^1}(1))$ are the sections vanishing at $\infty,
\,0\in\PP^1$ respectively. So its higher
cohomology has total dimension bounded by $O(k^n)$. But by pushing down
first to $X\times\PP^1$, then to $X$, this higher cohomology can be computed
as $\bigoplus_{j=0}^{ck}s^{ck-j}t^j\boxtimes H^{\ge1}_X(L^k\otimes\I_{\!Z}^{ck-j})$,
of total dimension $\sum_{j=1}^{ck}h^{\ge1}(L^k\otimes\I_{\!Z}^j)$.
\end{proof}

Fix $(\X,L)\to\C$ such that for $k\gg0$, $(\X,L^k)$ is a semi test configuration
for $(X,L^k)$. Given $\Z\subset\X$ a $\C^\times$-invariant
subscheme, denote its general fibre by $Z\subset X$ and central fibre $Z'\subset
\X_0$. Let $(\B_{Z'}(\X),\L_c)\rt{\pi}(\X,L)$
denote the blow up of $\X$ along $Z'$, with exceptional divisor $E$ and
line bundle $\L_c=\pi^*L-cE$. (As usual $c\in(0,\epsilon(Z'))$ or $c=\epsilon(Z')$
if $L^k\otimes\I_{\!Z'}^{\epsilon(Z')k}$ is saturated by
global sections for $k\gg0$.) Thus $\L_c$ is semi-ample by Proposition
\ref{ample}, making $(\B_{Z'}(\X),\L_c^k)$ a semi test configuration for
$(X,L^k)$ for $k\gg0$.

Then we have the following generalisation of Theorems \ref{degen} and \ref{ss}.

\begin{thm} \label{big}
In the above situation, suppose that the thickenings $j\Z\subset\X$ are flat
over $\C$ for all $j\in\mathbb N$, and the $\C^\times$-action on
$H^0_{\X_0}\!(L^k)$ has only weights
which lie between $-Ck$ and $0$, for some $C>0$. Then
$$
w\big(H^0_{\X_0}(L^k)\big)=w\big(H^0_\X(L^k)\big/tH^0_\X(L^k)\big)+O(k^n),
$$
$H^0_{(\B_{Z'}(\X))\_0}\!(\L_c^k)$ has only weights which lie between $-(C+c)k$
and $0$, and
$$
w\big(H^0_{(\B_{Z'}(\X))\_0}\!(\L_c^k)\big)=w\big(H^0_{\X_0}(L^k)\big)+w_k(Z)
+O(k^n).
$$
If $L$ is ample then the first $O(k^n)$ correction vanishes, and if in
addition either $c<\epsilon(Z')$ or the $\C^\times$-action on
$H^0_{\X_0}(L^k)$ is trivial then the second correction is $O(k^{n-1})$.
\end{thm}

\begin{proof}
$(\Y=\Proj\bigoplus_kH^0_\X(L^k),\O_\Y(1))$ is a polarised
family over $\C$ with general fibre $(\Proj\bigoplus_kH^0_X(L^k),\O(1))$.
It is flat because $H^0_\X(L^k)$ has no $t$-torsion, by the flatness of
$\X$. For $k\gg0$, by the flatness of $\Y$ and ampleness of
$\O_\Y(1)$, the central fibre has sections $H^0_{\Y_0}(\O(k))=
H^0_\Y(\O(k))\big/tH^0_\Y(\O(k))=H^0_\X(L^k)\big/tH^0_\X(L^k)$. Again by
flatness and ampleness, this has the same dimension
as for the general fibre, which is $h^0_X(L^k)$, which equals
$\chi\_X(L^k)+O(k^{n-1})$ by semi-ampleness. In turn by flatness and
semi-ampleness, this equals $\chi\_{\X_0}(L^k)+O(k^{n-1})=
h^0_{\X_0}(L^k)+O(k^{n-1})$. Therefore the inclusion
$$
\frac{H^0_\X(L^k)}{tH^0_\X(L^k)}\ \subseteq\ H^0_{\X_0}(L^k)
$$
has codimension $O(k^{n-1})$. This inclusion is $\C^\times$-equivariant,
and all weights on the right hand side lie between $-Ck$ and $0$ by assumption.
Therefore the total weights on the two vector spaces differ by at most
$O(k^n)$, as claimed. If $L$ is ample then by cohomology vanishing
$H^0_{\X_0}(L^k)=H^0_\X(L^k)\big/tH^0_\X(L^k)$ and the correction vanishes.
\\

To streamline notation we fix the convention in the proof of the second
result that
$ck+1=\infty$ so that, for instance, $\I_{\!Z'}^{ck}/\I_{\!Z'}^{ck+1}$ means
$\I_{\!Z'}^{ck}$ and we can deal with all of the terms in Proposition \ref{central}
uniformly. For $k\gg0$, this gives
\begin{equation} \label{HBl}
H^0_{(\B_{Z'}(\X))\_0}\!(\L_c^k)=H^0_{\X_0}\!\bigg(L^k\otimes
\frac{\I_{\!Z'\subset\X}^{ck}}{t\I_{\!Z'\subset\X}^{ck}}\bigg)=
\bigoplus_{i=0}^{ck}\ t^i\,H^0_{\X_0}\!\left(\!L^k\otimes
\frac{\I_{\!Z'}^{ck-i}}{\I_{\!Z'}^{ck-i+1}}\right).
\end{equation}
Since the weight on $t^i$ is $-i$ and lies between $-ck$ and $0$,  this
shows the weights on $H^0_{(\B_{Z'}(\X))\_0}\!(\L_c^k)$ indeed lie between
$-(C+c)k$ and $0$. Included in this $\C^\times$-module is
\begin{equation} \label{split2}
\bigoplus_{i=0}^{ck}\ t^i\frac{H^0_{\X_0}\!(L^k\otimes
\I_{\!Z'}^{ck-i})}{H^0_{\X_0}\!(L^k\otimes\I_{\!Z'}^{ck-i+1})}\,,
\end{equation}
of dimension $h^0_{\X_0}\!(L^k)$, which we have already noted above
is the same as $h^0_X(L^k)$ to $O(k^{n-1})$. The same working (applied
to $(\B_{Z'}(\X),\L_c)$ instead of $(\X,L)$) shows that
$h^0_{(\B_{Z'}(\X))\_0}\!(\L_c^k)$ also equals $h^0_X(L^k)+O(k^{n-1})$.
Thus the $-(C+c)k$-bound on weights means we can instead calculate the
weight on (\ref{split2}) at the expense of an $O(k^n)$ error.
If $L$ is ample and $c<\epsilon(Z')$ then $\L_c$ is also ample, so
$h^0_{(\B_{Z'}(\X))\_0}\!(\L_c^k)=h^0_X(L^k)=h^0_{\X_0}\!(L^k)$ and we
can calculate the weight on (\ref{split2}) without error.
Finally if $c=\epsilon(Z')$ and $L$ is ample
then as in Proposition \ref{error}, (\ref{HBl}) and (\ref{split2}) agree
for all but $i=0,\ldots eN$ (independent of $k$) so if the $\C^\times$-action
on $H^0_{\X_0}(L^k)$ is trivial then their weights differ by
$\le\sum_{i=0}^{eN}ih^1(L^k\otimes\I_{\!Z'}^{ck-i+1})\le eN.O(k^{n-1})=O(k^{n-1})$.

Define $V^0:=H^0_{\X_0}\!(L^k)$, and
\begin{equation} \label{filt}
V^p:=H^0_{\X_0}\!(L^k\otimes\I_{\!Z'}^p)\subseteq V^{p-1}\subseteq\ldots
\subseteq V^0.
\end{equation}
Let $V^0=\bigoplus_{j=-Ck}^{\,0}V^{0,j}$ be its weight space decomposition.
Given such a splitting $\bigoplus_jV^{0,j}$ of a vector space
$V^0$, the generic subspace $V^p\subset V^0$ is not generated
by the pieces $V^{p,j}:=V^p\cap V^{0,j}$; i.e.\ $V^p\supsetneq
\bigoplus_jV^{p,j}$. But if $V^0$ has a $\C^\times$-action whose weight space
decomposition is $\bigoplus_jV^{0,j}$, and each $V^p\subset V^0$ is
$\C^\times$-invariant, then indeed
\begin{equation} \label{filt1}
V^p=\bigoplus_jt^jV^{p,j} \qquad\text{and}\qquad \frac{V^p}{V^{p+1}}\cong
\bigoplus_jt^j\frac{V^{p,j}}{V^{p+1,j}}\,.
\end{equation}
This holds here since $\I_{\!Z'}^p$ is $\C^\times$-invariant.

So the splitting $\bigoplus_{i=0}^{ck}t^iV^{k-i}\big/V^{k-i+1}$ (\ref{split2})
further splits as
\begin{equation} \label{split9}
\bigoplus_{i,j}t^i\frac{V^{ck-i,j}}{V^{ck-i+1,j}}\,.
\end{equation}
Defining $h^{a,b}:=\dim V^{ck-a,b}/V^{ck-a+1,b}$, the weight on the determinant
of (\ref{split9}) is $\sum_{i,j}(-i+j)h^{i,j}$.

By (\ref{filt1}), $\sum_jh^{i,j}=\dim V^{ck-i}/\dim V^{ck-i+1}$, which by (\ref{filt})
is $h^0_{\X_0}\!(L^k\otimes\I_{\!Z'}^{ck-i})-h^0_{\X_0}\!(L^k\otimes\I_{\!Z'}^{ck-i+1})$.
So as in (\ref{noco}), $-\sum_{i=0}^{ck}i\sum_jh^{i,j}=\sum_{i=1}^{ck}h^0_{\X_0}
(L^k\otimes\I_{\!Z'}^i)\,-ckh^0_{\X_0}\!(L^k)$.

Similarly, the fact that $V^{i+1,j}\subset V^{i,j}\subset\ldots\subset V^{0,j}$
filters $V^{0,j}$ means that $\sum_ih^{i,j}=\dim V^{0,j}$. Therefore
$$
w\big(H^0_{(\B_{Z'}(\X))\_0}\!(\L_c^k)\big)=\sum_{i,j}(-i+j)h^{i,j}=
-\sum_{i=0}^{ck}i\bigg(\!\sum_jh^{i,j}\bigg)
+\!\sum_{j=-Ck}^0j\bigg(\!\sum_ih^{i,j}\bigg)+O(k^n)
$$
equals
$$
\sum_{i=1}^{ck}h^0_{\X_0}\!(L^k\otimes\I_{\!Z'}^{ck-i})-ckh^0_{\X_0}\!(L^k)\
+\sum_{j=-Ck}^0j\dim V^{0,j}+O(k^n).
$$
We can replace $h^0$ by $\chi$ at the expense of $O(k^n)$ by Lemma \ref{wkerror}
(and at the expense of $O(k^{n-1})$, by Proposition \ref{error}, if $L$
is ample). But $\I_{\!Z'}^j\subset\O_{\X_0}$ sits in a flat family with
$\I_{\!Z}^j\subset\O_X$, in which Euler characteristics are preserved, so the above
weight is
\begin{multline*}
\sum_{j=1}^{ck}\chi\_X(L^k\otimes\I_{\!Z}^j)-ckh^0_X(L^k)+w\big(
H^0_{\X_0}(L^k)\big)+O(k^n) \\
=w_k(Z)+w\big(H^0_{\X_0}(L^k)\big)+O(k^n),
\end{multline*}
by (\ref{wkchi}). And if $L$ is ample and if either $c<\epsilon(Z')$ or the
$\C^\times$-action on $H^0_{\X_0}(L^k)$ is trivial then the correction
is $O(k^{n-1})$.
\end{proof}

\section{Towards a converse} \label{converse}

To apply Theorem \ref{big} to Corollary \ref{ref} to express weights of test
configurations (\ref{any}) as sums of $w_k(Z)$s requires flatness of the
thickenings $j(\overline{Z_i\times\C})$ in $\X^i$. To help achieve this by
applying
Proposition \ref{flat} we assume that $X$ is reduced and pass to a blow up
$p\colon\widehat X\to X$ on which the (pullbacks of the) $Z_i$ are Cartier
divisors $D_i$: $p^*\I_{Z_i}=\O(-D_i)$. In fact
by Hironaka's resolution of singularities we may take the $D_i$ to have simple
normal crossing (snc) support. That
is there are smooth reduced divisors $\{F_j\}\subset\widehat X$ such that
$F=\cup_jF_j$ has simple normal crossings and $D_i=\sum_jm_{ij}F_j$ for each
$i$ and some nonnegative integers $m_{ij}$. We pull $L$ back
to $\widehat X$, and construct the families $(\XX^s,L_s)\to\C$
as before, using $D_s\subset\widehat X$ in place
of $Z_s\subset X$. There are equivariant surjective maps $\XX^s\to\X^s$
which identify (by pullback) the $L_s$
line bundles, defined by pulling back sections of $L^k\otimes\I_i^j$ on $X$
to sections of $L^k\otimes\O(-jm_iD_i)$ on $\widehat X$.

\begin{thm} \label{bigger}
Suppose that $X$ is normal, and fix an arbitrary test configuration $(\Y,\O_\Y(1))$
(\ref{any}) with associated subschemes $Z_0\subseteq\ldots\subseteq Z_{r-1}\subset
X$ (\ref{ideals}). Suppose that there is a resolution $D_0\subseteq\ldots\subseteq
D_{r-1}\subset\widehat X$ with divisors $D_i$ satisfying condition
(\ref{reduced}). (For instance if the connected components of the snc divisors
$D_i$ have the same multiplicities locally, i.e.\ $m_{ij}\in\{0,m\}$ for some
locally constant $m$ and all $i$. E.g. if the $D_i$ are all reduced, which
is the $m=1$ case.) Then
$$
w(H^0_{\Y_0}\!(\O_\Y(k)))=\sum_{i=0}^{r-1}w_k(Z_i)+O(k^n).
$$ 
\end{thm}

\begin{proof}
$(\XX^r,L_r)$ dominates $(\X^r,L_r)$ which in turn dominates
$(\B_{I_r}(X\times\C),\L_1)$ and so $(\Y,\O_\Y(1))$,
and the space of sections of the $k$th power of the line bundle on the
general fibre is $h^0_X(L^k)$ for all three, by the normality of $X$.
Therefore by Proposition \ref{stein}
and the first part of Theorem \ref{big} we may calculate $w(H^0_{(\widehat
\X^r)\_0}\!(L_r^k))$ at the expense of an $O(k^n)$ error.

We apply Theorem \ref{big} inductively to $\X=\XX^i$ with $c=1$ (which is
$\le\epsilon(Z_i)=\epsilon(D_i)=\epsilon(D^{(i)}_i)$ by Corollary
\ref{seshadriIr}) and $\Z=\overline{D_i\times\C}$ with central fibre
$D_i^{(i)}$. The condition (\ref{reduced}) on the $D_i$ guarantees that each
$j(\overline{D_i\times\C})$ is flat over $\C$ by Proposition \ref{flat}.

The induction starts with $\widehat X\times\C$, for which all weights are
zero so trivially satisfy the $-Ck$ bound. Theorem \ref{big} then ensures
the induction continues to compute the weight as
$$
\sum_{i=0}^{r-1}w_k(D_i)+O(k^n).
$$
Since $X$ is normal, all sections of $p^*L$ on $\widehat X$ are pullbacks
from $X:\ H^0_{\widehat X}(L^k)\cong H^0_X(L^k)$. Thus the same is true of
those sections vanishing on the pullback of $Z_i:\ H^0_{\widehat X}
(L^k\otimes\I_{\!D_i}^j)\cong H^0_X(L^k\otimes\I_{\!Z_i}^j)$. Thus $w_k(D_i)
=w_k(Z_i)$ by their definition (\ref{wkZ}), as required.
\end{proof}

In particular we can now calculate the leading order term of the weight if
the snc divisors
$D_i$ in a resolution of singularities $\widehat X$ of $(X,Z_i)$ are
reduced. The next case we consider is when they have multiplicities which
can vary with $i$ but are still locally constant over their snc support.
Here the relevant flatness result does not hold, but we will find that it
does after performing a basechange and normalisation.

So we consider the case when $I_r$ (\ref{ideals}) is locally of the form
\begin{equation} \label{ideal}
\I_D^{p_1}+t\I_D^{p_2}+\ldots+(t^r),
\end{equation}
for some reduced snc divisor $D$. These $p_i$ may
vary over the different connected components of $D$, but since the total
weight is a sum over contributions from each connected component we can
calculate the weight of each separately.

Pick a local function $z$ generating the ideal $\I_D$, so that (\ref{ideal})
is $(z^{p_1})+t(z^{p_2})+\ldots+(t^r)$. In the concave hull of the points
$(p_i,i),\,i=1,\ldots,r$ in the $(z,t)$-plane, we choose extremal vertices
$(k_i,\rho_i),\,i=1,\ldots,l$ so that they form a
concave set with the same concave hull (and $(k_1,\rho_1)=(p_1,0),\ (k_l,\rho_l)
=(0,r)$). This defines a new ideal
\begin{equation} \label{ideal2}
(z^{k_1})+t^{\rho\_2}(z^{k_2})+\ldots+(t^{\rho\_l}),
\end{equation}
with the same integral closure as (\ref{ideal}) since in this situation
taking integral closures corresponds to taking concave hulls, and an ideal
saturates its integral closure. In the next theorem we decorate $w_k(D)$
(\ref{wkZ}) as $w_k(D,c)$ with the value $c$ that determines the line bundle
$\L_c=L-cE$ on the blow up.

\begin{thm} \label{divs}
Set $m_i=\frac{\rho_{i+1}-\rho_i}{k_i-k_{i+1}},\ i=1,\ldots l$, and $m_0=0$.
If $L$ is ample, then the total weight of the blow up in the ideal (\ref{ideal})
is the sum over the connected components $D$ of
$$
-\sum_{i=1}^{l-1}(m_i-m_{i-1})w_k(D,k_i)-ak^n+O(k^{n-1}),
$$
for some $a\ge0$. If $L$ is semi-ample, the above expression is correct
to $O(k^n)$.
\end{thm}

\begin{proof}
We start by proving the weaker estimate for $L$ semi-ample.
First blow up $X\times\C$ in $D\times\{0\}$ (i.e.\ locally in the ideal $(z)+(t)$),
then in $D'$ (recall this is
the central fibre of $\overline{D\times\C}$ in the blow up), then $D''=D^{(2)}$,
etc. up to $D^{(j-1)}$. As in Theorem \ref{pr} we denote by $E_i$ the
pullback of the exceptional divisor of the $i$th blow up, and by $s_i$
the canonical section of $\O(E_i)$ vanishing on $E_i$. Then we claim
that the pushdown of $\O(-p_1E_1-\ldots-p_jE_j),\ p_1\ge p_2\ge\ldots\ge
p_j$, to $X\times\C$ is the ideal
\begin{multline} \label{induction}
(z^{p_1})+\ldots+t^{p_1+p_2-2p_2}(z^{p_2})+\ldots+t^{p_1+p_2+p_3-3p_3}(z^{p_3})
\\ +\ldots+\qquad\cdots\qquad+\ldots+
t^{p_1+\ldots+p_j-jp_j}(z^{p_j})+(t^{p_1+\ldots+p_j}).
\end{multline}
Here ``$\ldots$" means ``all convex combinations in between", i.e.\ the integral
closure of the ideal generated by the named terms.
(So $(f)+\ldots+(g)$ includes all monomials $h$ such that there exist $\lambda,\mu
\in\mathbb N$ with $h^{\lambda+\mu}=f^\lambda g^\mu$.) That is, we claim
the sections of $\O(-p_1E_1-\ldots-p_jE_j)$ over $X\times\C$ are $s_1^{-p_1}\ldots
s_j^{-p_j}$ times by the above ideal. This is standard but
fiddly to prove, and is best done by Newton diagram. We give
an unenlightening proof; the reader is advised to skip straight to the
Newton diagram (Figure \ref{diag}) for the special case below.

We prove (\ref{induction}) inductively alongside the claim that the complement
of the divisor $t/s_j=0$ (which is the proper transform, in the $j$th blow
up, of the central fibre of the $(j-1)$th blow up) is affine over $X\times\C$
with coordinate ring
\begin{equation} \label{affcoord}
\O\_{\!X\times\C}\left[\frac z{t^j}\right].
\end{equation}
(In particular, at a smooth point of $D$ where $z$ is a local coordinate
in an analytic coordinate system $(y_1,\ldots y_{n-1},z)$ for $X$, we see
that this part of the $j$th blow up is $\Spec\C[y_1,\ldots,y_{n-1},z,t,z/t^j]
=\Spec\C[y_1,\ldots,y_{n-1},Z,t]$, where $Z=z/t^j$, and so is locally
isomorphic to $X\times\C$ with the proper transform of $(z^k=0)$ being $(Z^k=0)$.)

For the first blow up, the sections of $O(-kE_1)$ are $(z,t)^ks_1^{-k}$ (there
are no more because the exceptional divisor is a $\PP^1$-bundle over $D$).
This is in agreement with (\ref{induction}),
and the proper transform of the central fibre is $t/s_1=0$. On its complement,
$t^k/s_1^k$ trivialises $O(-kE_1)$; dividing by it identifies the sections
$((z)+(t))^ks_1^{-k}$ with $\O\_{\!X\times\C}.((z/t)^k+(z/t)^{k-1}+\ldots+(z/t)+1)$.
Taking the limit as $k\to\infty$ gives the coordinate ring $\O\_{\!X\times\C}
\left[\frac zt\right]$ claimed (\ref{affcoord}). So the induction starts
at $j=1$.

At the $j$th stage the coordinate ring (\ref{affcoord}) shows that $D^{(j)}$
has ideal $(Z,t)$, where $Z=z/t^j$. Let $\pi$ denote the $(j+1)$th blow
up. By the first step of the induction, $\{t/s_{j+1}\ne0\}$ has coordinate
ring augmented by $Z/t=z/t^{j+1}$;
i.e.\ by induction $\O\_{\!X\times\C}\left[\frac z{t^j}\right]\left[\frac
z{t^{j+1}}\right]=\O\_{\!X\times\C}\left[\frac z{t^{j+1}}\right]$.

$\pi_*\O(-p\_{j+1}E_{j+1})$ is the ideal $((Z)+(t))^{p\_{j+1}}=
t^{-jp\_{j+1}}(z^{p\_{j+1}})+\ldots+(t^{p\_{j+1}})$; i.e.\ the sections of
$\O(-p\_{j+1}E_{j+1})$ over $\X^j$ are $s_{j+1}^{-p\_{j+1}}$ times by sections
of this ideal. Multiplying by the trivialising section 
$t^{p_1}s_1^{-p_1}\ldots t^{p_j}s_j^{-p_j}$ of $\O(-p_1E_1-\ldots-p_jE_j)$
shows that over our affine piece,
$$
\pi_*\O(-p_1E_1-\ldots-p\_{j+1}E_{j+1})=\big(\pi_*\O(-p\_{j+1}E_{j+1})\big)\otimes
\O(-p_1E_1-\ldots-p_jE_j)
$$
is the ideal
$$
\big[t^{p_1+\ldots+p_j-jp\_{j+1}}(z^{p\_{j+1}})+\ldots+(t^{p_1+\ldots+p_j+p\_{j+1}})
\big].s_1^{-p_1}\!\ldots s_j^{-p_j}\subseteq\O(-p_1E_1-\ldots-p_jE_j).
$$
Thus the sections of $\O(-p_1E_1-\ldots-p\_{j+1}E_{j+1})$ are the sections
of $\O(-p_1E_1-\ldots-p_jE_j)$ which lie in the above ideal, i.e.\ the intersection
of (\ref{induction}) with the ideal $t^{p_1+\ldots+p_j-jp\_{j+1}}(z^{p\_{j+1}})
+\ldots+(t^{p_1+\ldots+p_j+p\_{j+1}})$. But this is (\ref{induction}) with
$j$ replaced by $j+1$, completing the induction.

We first assume that the $m_i$ are all integers. Then the ideal (\ref{ideal2})
has integral closure of the form (\ref{induction}), on taking $p_1,\ldots,p_{m_1}$
all equal
\begin{figure}[h]
\center{\input{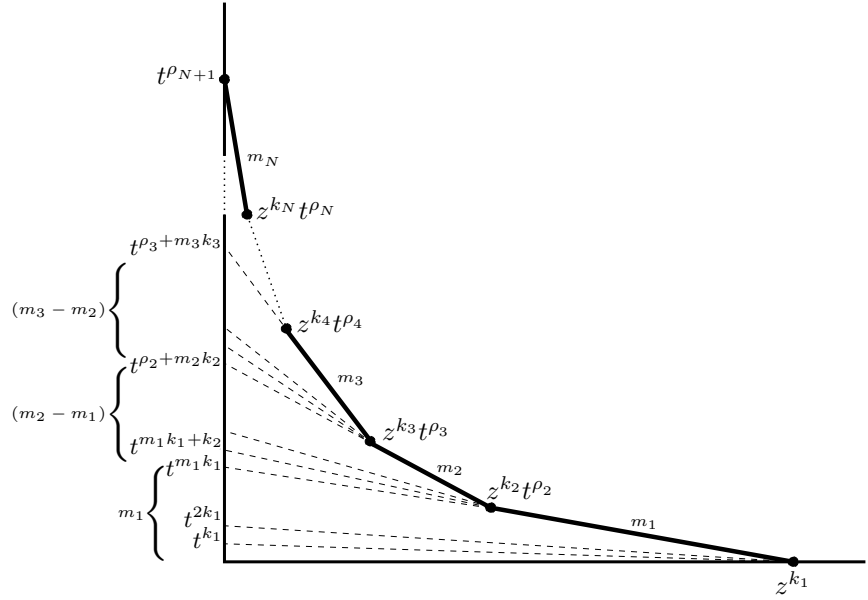} \caption{Newton diagram for the blow up in
the ideal (\ref{ideal2}) \label{diag}}}
\end{figure}
to $k_1$, then $p_{m_1+1},\ldots,p_{m_2}$ all equal to $k_2$, and so
on, up to $p_{m\_{N-1}+1},\ldots,p_{m\_N}$ all equal to $k_N$.
This is illustrated in Figure \ref{diag}, the Newton diagram of the $(z,t)$
plane, with the $m_i$s being (minus) the gradients of the bold lines.
Taking the integral closure, i.e.\ including the monomials ``$\ldots$" in
(\ref{induction}), corresponds to including all integral points
both on and above the line to lie in the ideal. Replacing $z$ by $Z=z/t^j$
and multiplying by the trivialising section $t^{p_1+\ldots+p_j}$ in the
above working corresponds to the integral affine transformation that locally
takes one corner of the bold line into the $(z,t)$-axes.

So we may calculate the weight of this sequence of blow ups, $m_1$ times
in the central fibre of $\overline{D\times\C}$ with weight $c=k_1$, then
$(m_2-m_1)$ times in the central fibre of $\overline{D\times\C}$ with weight
$c=k_2$, etc. (The weight here just means the coefficient of the exceptional
divisor in the line bundle $L-cE$ we use, as usual.) The flatness criterion
(\ref{reduced}) is trivially satisfied so that by Theorem \ref{big} we
may calculate the weight to be that claimed. This differs from the weight
of the blow up in (\ref{ideal}) by a $-ak^n+O(k^{n-1})$ correction by Proposition
\ref{stein} since the blow up in the integral closure (\ref{ideal2})
is the normalisation of the blow up in (\ref{ideal}).

If $L$ is ample then we can improve the estimate. In the first blow up,
$\X=X\times\C$ is the trivial product configuration, so the
$\C^\times$-action on $H^0_{\X_0}(L^k)$ is trivial and we can use the better
estimate of Theorem \ref{big}. Next we group the first $m_1$ blow ups together
as one blow up in $((z)+(t^{m_1}))^{k_1}$; this has the advantage
that we do no blowing down (in fact the previous $m_1$ blow ups blow down
to this). Since this blow up is the $t\mapsto t^{m_1}$ basechange of the
blow up in the ideal $((z)+(t))$ (with $c=k_1$), it has weight $m_1w_k(D,k_1)$,
the same as the sum of the $m_1$ blow ups we performed above. Similarly we
group the next
$m_2$ blow ups together, blowing up with $c=k_2$ in the ideal generated
by $t^{m_2}$ and the ideal of $\overline{D\times\C}$, and calculate its
weight as the $t\mapsto t^{m_2}$ basechange of a blow up we already know.
Inductively we end up with the same formula for the weights, but with
the added $O(k^{n-1})$ accuracy of Theorem \ref{big} coming from the fact
that (after the first blow up) there are no blow downs, so $c$ is less
than the Seshadri constant of the relevant $D^{(i)}$ at each stage.

Finally, if the $m_i$ are not integers, we simple replace $t$ by $t^M$
in (\ref{ideal}), i.e.\ we basechange our test configuration, where $M$
clears the denominators of all the $m_i$. This replaces all the $m_i$ by
the integers $Mm_i$ while multiplying the $\C^\times$-weight by $M$.
Substituting the $Mm_i$ into our formula for the weights gives the weight
of this new test configuration with a $-ak^n+O(k^{n-1})$ correction (coming
from taking
the integral closure of this new ideal; replacing the test configuration
with its normalisation and using Proposition \ref{stein}); dividing by $M$
gives the weight of the original test configuration for any $m_i$.
\end{proof}

This just leaves the case of where the divisors $F_j$ in the snc divisors
$D_i$ in a resolution $\widehat X$ intersect with differing multiplicities
(the simplest example
being $\I_{D_0}=(x^2y)$ and $\I_{D_1}=(x)$ locally). This of course cannot
happen for curves, so we have

\begin{cor} \label{Kslopecurve}
A smooth curve $(X,L)$ is K-(semi/poly)stable if and only if it is slope
(semi/poly)stable.
\end{cor}

\begin{proof} 
Smooth curves are normal, and for the resolution of singularities $\widehat
X$ we of course take $X$ itself, so we can apply the stronger form of Theorem
\ref{divs}. Since the Donaldson-Futaki invariant only uses the coefficients
$b_0$ and $b_1$ of $w(k)=b_0k^{n+1}+b_1k^n+O(k^{n-1})$, this implies that
the Futaki
invariant of an arbitrary test configuration (\ref{any}) is $\ge$ a positive
linear combination
of Futaki invariants of (the deformation to the normal cone of) subschemes.
Slope stability implies that these are all positive.
\end{proof}

This result makes it trivial to understand K-stability for smooth curves;
see Theorem \ref{smoothcurvesareKstable}.

\section{Chow stability} \label{sec:chow}

Mumford's notion of Chow (semi)stability \cite{Mu} of $(X,\O_X(1))\subseteq\PP^N$,
for \emph{fixed} $N$, is the simplest form of stability to calculate (as
opposed to \emph{asymptotic} Chow stability (\ref{chow}), which is second
only to asymptotic Hilbert stability in difficulty). It is also
useful in algebro-geometric applications since it is a genuine GIT
notion giving projective (and so proper and separated) moduli spaces.

For any subscheme $Z\subseteq X$, we define the polynomial
$a_0(x)$ by $\chi\_X(\I_{\!Z}^{xk}(k))=a_0(x)k^n+a_1(x)k^{n-1}+\ldots$ for
$k\gg x^{-1}>0$ as before. We also define $a_0$ by $h^0(\O_X(k))=a_0k^n+
a_1k^{n-1}+\ldots$\,; by (\ref{a0}), $a_0=a_0(0)$.

For any subscheme $Z\subseteq X\subseteq\PP^N$ and \emph{integer}
$0<c\le\epsilon(Z)$, we define the \emph{Chow slope} of $\I_Z$ to be
$$
Ch_c(\I_Z):=\frac{\sum_{i=1}^ch^0_{\PP^N}(\I_{\!Z}^i(1))}{\int_0^ca_0(x)dx}\,,\qquad
Ch(\I_Z)=\max_{\mathbb N\ni c\le\epsilon(Z)}\big(Ch_c(\I_Z)\big)\ \in\,[-\infty,\infty].
$$
(Here we define max of the empty set to be $-\infty$, and division by 0
to give $\infty$.) Setting $Z=\emptyset$ gives $Ch(X)=Ch(\O_X)$ as
$$
Ch(X):=\frac{h^0_{\PP^N}(\O(1))}{a_0}=\frac{N+1}{a_0}\,.
$$
If $X\into\PP^N$ is the Kodaira embedding of $X$ in $\PP(H^0(\O_X(1))^*)$,
i.e.\ if $H^0_{\PP^N}(\O(1))\to H^0_X(\O(1))$ is an isomorphism, then
$h^0_{\PP^N}(\I_Z(1))=h^0_X(\I_Z(1))$ for all $Z\subseteq X$. This can
be arranged, for fixed $(X,L)$, by taking a sufficiently large multiple
$L^k=:\O_X(1)$ of the polarisation and setting $\PP^N=\PP(H^0_X(L^k)^*)$.
In this situation the above slopes can be written intrinsically in terms
of $(X,\O_X(1))$ as
$$
Ch_c(\I_Z)=\frac{\sum_{i=1}^ch^0_X(\I^i_{\!Z}(1))}{\int_0^ca_0(x)dx}
\qquad\text{and}\qquad Ch(X)=\frac{h^0_X(\O(1))}{a_0}\,.
$$

We say that $X\subset\PP^N$ is \emph{Chow slope stable} if $Ch(\I_Z)<Ch(X)$
for all nonempty $Z\subseteq X$, and Chow slope semistable if $Ch(\I_Z)\le
Ch(X)$. We have
\begin{equation} \label{sss}
Ch_c(\I_Z)<Ch(X) \quad\Longleftrightarrow\quad
Ch(X)<Ch_c(\O_Z):=\frac{\sum_{i=1}^c\big(N+1-h^0(\I_{\!Z}^i(1))\big)}
{\int_0^c\!\tilde a_0(x)dx}\,.
\end{equation}

\begin{thm} \label{thmchow}
If $X\subseteq\PP^N$ is Chow (semi)stable then it is Chow slope (semi)stable.
\end{thm}

\begin{proof}
Choose a basis of $H^0_{\PP^N}(\O(1))$ compatible with the filtration
$H^0_{\PP^N}(\I_{\!Z}^c(1))\subseteq H^0_{\PP^N}(\I_{\!Z}^{c-1}(1))\subseteq\ldots
\subseteq H^0_{\PP^N}(\I_Z(1))\subseteq H^0_{\PP^N}(\O(1))$, so that the
first $p_i:=h^0_{\PP^N}(\I_{\!Z}^i(1))$ elements are contained in
$H^0_{\PP^N}(\I_{\!Z}^i(1))$. The corresponding hyperplanes $H_1,H_2,\ldots,$
$H_{N+1}$ and subschemes $Z_i:=
X\cap H_1\cap\ldots\cap H_i$ therefore satisfy $Z_{p_i}\supseteq iZ$ (i.e.\
$\I_{Z_{p_i}}\subseteq\I_{\!Z}^i$).

Choose weights $\rho_1=0=\ldots=\rho_{p_c},\ \rho\_{p\_
c+1}=1=\ldots=\rho_{p_{c-1}},\
\ldots,\ \rho\_{p\_1+1}=c=\ldots=\rho_{p_0}$ so that ideal on $X\times\C$,
\begin{equation} \label{IN}
t^{\rho_1}\I_{Z_1}+t^{\rho_2}\I_{Z_2}+\ldots+t^{\rho\_{\!N}}\I_{Z_N}
+(t^{\rho\_{\!N+1}})
\end{equation}
equals
\begin{equation} \label{wobbly}
I=\I_{Z_{p_c}}+t\I_{Z_{p_{c-1}}}+\ldots+t^{c-1}\I_{Z_{p_1}}+(t^c),
\end{equation}
which by construction is contained in the ideal
\begin{equation} \label{straight}
\I_{\!Z}^c+t\I_{\!Z}^{c-1}+\ldots+t^{c-1}\I_Z+(t^c)=(\I_Z+(t))^c.
\end{equation}

$X\subset\PP^N$ is Chow stable for the $\C^\times$-action which has
weight $\rho_i$ on the $i$th vector of our basis of $H^0_{\PP^N}(\O(1))$.
Mumford (\cite{Mu} Theorem 2.9) shows that this is equivalent to the inequality
\begin{equation} \label{mummy}
-a<\frac{a_0}{N+1}\sum_{i=1}^{N+1}\rho_i,
\end{equation}
where $a$ is the Chow weight of the blow up of $X\times\C$ in
the ideal $I$ (\ref{wobbly}), i.e.\ $ak^{n+1}+O(k^n)$ is the total weight
of the induced $\C^\times$-action on $H^0_{(\B_I(X\times\C))\_0}\!(L^k(-kE))$.
(This is only one half of Mumford's result, the
harder part being his computation of $a$ in terms of
ideals on $X\times\C$, which we do not use. The reader wishing
to compare our conventions with Mumford's should rewrite the above as
$-(n+1)!\,a<(n!\,a_0)\frac{n+1}{N+1}\sum\rho_i$,
replace $n$ with $r$, $N$ with $n$, and $\rho_i$ with $\rho\_{N-i}$.
Finally the sign arises from our convention (\ref{convention}) that if
$g$ acts on $V$ then the induced action on its functions $S^KV^*$ is by
$(S^Kg^*)^{-1}$; Mumford calculates the weight of $S^Kg^*$.) Notice this
is Theorem \ref{thm:instability} applied to (\ref{normalised}) on setting
$r=1,\,h^0(L^r)=N+1$ and $w(1)=\sum_i\rho_i$. 

Since $c\ge\epsilon(Z)$, we know the Chow weight on the blow up in $(\I_Z+(t))^c$
(Proposition \ref{integrals}). But $I$ is
contained in $(\I_Z+(t))^c$ (\ref{straight}), so the weight on the latter
is more negative (for instance Mumford's formula (\cite{Mu} Theorem 2.9)
for the Chow weight shows this), giving the inequality
$$
a\le\int_0^ca_0(x)dx-ca_0.
$$
Thus (\ref{mummy}) gives
$$
\frac{N+1}{a_0}\!\left(\!ca_0-\int_0^ca_0(x)dx\!\right)\!<
\sum_{i=1}^{N+1}\rho_i=1(p_{c-1}-p_c)+2(p_{c-2}-p_{c-1})+\ldots+c(p_0-p_1)
\vspace{-10pt} $$$$ \hspace{3cm}
=-p_c-p_{c-1}-\ldots-p_1+cp_0=-\sum_{i=1}^ch^0_{\PP^N}(\I_{\!Z}^i(1))+c(N+1),
$$
which is
\begin{equation} \label{r1}
\frac{\sum_{i=1}^ch^0_{\PP^N}(\I^i_{\!Z}(1))}{\int_0^ca_0(x)dx}
<\frac{N+1}{a_0}\,.
\end{equation}
Semistability is similar, replacing (\ref{mummy}) by the non strict inequality.
\end{proof}

\begin{rmks}
Setting $c=1$ gives the inequality $\frac{N+1}{a_0}\int_0^1a_0(x)dx>h^0(\I_Z(1))$,
and replacing $Z$ by $kZ$ then gives $\frac{N+1}{a_0}\frac1k\int_0^ka_0(x)dx>
h^0(\I_{\!Z}^k(1))$. We could have gotten this alternative slope-type inequality
directly from stability for the $\C^\times$-action that had all of the above
$\rho_i$ equal to zero (for $i\le p_k$) or one ($i>p_k$). However, it is
strictly weaker than our slope inequality $\frac{N+1}{a_0}\int_0^ka_0(x)dx>
\sum_{i=1}^kh^0(\I_{\!Z}^i(1))$ (\ref{thmchow}) since this last sum is clearly
$\ge kh^0(\I_{\!Z}^k(1))$. \smallskip

For $X$ semistable, taking $Z=X,\,c=1$ (so $a_0(x)\equiv0$) shows that
$h^0_{\PP^n}(\I_X(1))=0$, i.e.\ $X$ is not contained in any hyperplane and
$H^0_{\PP^N}(\O(1))$ injects into $H^0_X(\O(1))$. To make this injection
an isomorphism requires the assumption that $X\subseteq\PP^N$ is a Kodaira
embedding. \smallskip

Notice that (\ref{r1}) is more-or-less the $k^{n+1}$
coefficient of (\ref{weight}) (with $r=1$ and rearranged). There can be
higher cohomology corrections since $I$ (\ref{wobbly}) is not the same as
$(\I_Z+(t))^c$ (\ref{straight}); these disappear for $r\gg0$.
\end{rmks}

We now want to understand to what extent slope Chow stability should
imply Chow stability. This involves
demonstrating the inequality (\ref{mummy}) for \emph{all} linearly independent
sequences of hyperplanes $H_1,H_2,\ldots$ and \emph{all} choices of weights
$0=\rho_1\le\rho_2\le\ldots\le\rho\_{N+1}$. As before we set $Z_i=X\cap
H_1\cap\ldots H_i,\ \I_i:=\I_{Z_i}$.

The idea is to relate the weight of the associated $\C^\times$-action (with
weight $\rho_i$ on the $i$th vector of our basis of $H^0_{\PP^N}(\O(1))$)
to a sum of weights of $\C^\times$-actions of the standard form considered
in Theorem \ref{thmchow}. The problem is that we have noted
that we can only express weights as a sum of weights of the deformation
to the normal cone of subschemes if either
condition (\ref{reduced}) holds or the multiplicities of $\I_i$ are locally
constant (at least on some resolution where they are ncds). But in either
of these good cases we can demonstrate the general procedure of passing
from slope stability to stability.

So let us assume for illustration that each $\I_i$ has the local form $O(-k_iD)$
for some reduced Cartier divisor $D$ and number $k_i$ constant on $D$.
(We have seen, as in Theorem \ref{bigger}, one can also pass to a resolution
$\widehat X$ if necessary to make the $Z_i$ sncds to calculate Chow weights
(the $k^{n+1}$-coefficient of Hilbert weights), so the assumption is not
so restrictive, and is sufficient to deal with curves.)
That is, $I=t^{\rho_1}\I_1+t^{\rho_2}\I_2+\ldots+
(t^{\rho\_{\!N+1}})$ has the local form
\begin{equation} \label{deali}
I=\I_D^{k_1}+t^{\rho_2}\I_D^{k_2}+\ldots+(t^{\rho\_{l}}),
\end{equation}
with $0=\rho_1\le\rho_2\le\ldots\le\rho\_{\!N+1}$,
and $k_1\ge k_2\ge\ldots>k_l=0=k_{l+1}=\ldots=k_{N+1}$ (so $l\le N+1$ is
the smallest number with $k_l=0$). These $k_i$ and $l$ may vary as $D$ ranges
over the connected components of $Z_1$.

\begin{thm} \label{uniform}
If $(X,\O(1))$ is slope Chow stable then the weight $a$ (\ref{mummy})
satisfies $-\frac{N+1}{a_0}\,a<\sum_{i=1}^{N+1}\rho_i+\delta$,
where $\delta=\sum_{i=1}^N(\rho_{i+1}-\rho_i)h^1(\I_i(1))$.
\end{thm}

This is the inequality  $-\frac{N+1}{a_0}\,a<
\sum_{i=1}^{N+1}\rho_i$ required for Chow stability (\ref{mummy}), modulo
some $h^1$ corrections (which we estimated away in the K-stability
analogue). In (\ref{uniformcurve}) the result is strengthened slightly and
the correction estimated on curves to prove their asymptotic
Chow stability.

\begin{proof}
The $k^{n+1}$-coefficient $a$ of
the weight of our $\C^\times$-action is the same as that on the blow
up of $\widehat X\times\C$ in $I$ (\ref{deali}). In Theorem \ref{divs} we
calculated this to be a sum $a=\sum_Da_D$ over
the connected components $D$ of its support, where
\begin{equation} \label{a}
a_D=-\sum_{i=1}^{l-1}{}'\left(\!\frac{\tilde\rho_{i+1}-\tilde\rho_i}{k_i-k_{i+1}}-
\frac{\tilde\rho_i-\tilde\rho_{i-1}}{k_{i-1}-k_i}\right)\!
\int_0^{k_i}\!\tilde a_0^D(x)dx.
\end{equation}
Here $\tilde\rho_i$ is defined uniquely by requiring that $(k_i,\tilde\rho_i)$
lies on the boundary of the concave hull of the set of points $(k_i,\rho_i)_{i=1}^l$
in the $(k,\rho)$-plane. Thus the $\tilde\rho_i$ need not be integers but,
for instance, $\tilde\rho_1=\rho_1$ and $\tilde\rho_i=\rho_l$ for all $i\ge
l$. More generally, $\tilde\rho_i\le\rho_i$ for all $i$. Concavity of the
$(k_i,\tilde\rho_i)$ ensures that any term in the above sum with a zero
in the denominator also has zero in the numerator; the prime $'$ on the summation
sign signifies that we ignore these $\frac00$ terms (the terms with $k_i=k_{i+1}$)
in the sum; equivalently
we set $\frac00:=0$. In the first term we set $\tilde\rho_0:=0$.

The Seshadri constant of $D$ is $\ge k_i$ for
all $i$ since $\I_D^{k_i}$ is locally the intersection of
a sequence of hyperplanes, so $\I_D^{k_i}(1)$ is globally generated
near $D$, so $\I_D^{rk_i}(r)$ is too. Therefore Chow slope stability
for $D\subset X$ (\ref{sss}) gives the inequalities
$$
\frac{N+1}{a_0}\int_0^{k_i}\!\tilde a_0^D(x)dx<\sum_{i=1}^{k_i}\big(
N+1-h^0(\I^i_D(1))\big).
$$
Since all of the integrals in (\ref{a}) have coefficients which are $\ge0$
by the concavity of the $(k_i,\tilde\rho_i)$, we obtain
\begin{eqnarray} \nonumber
-\frac{N+1}{a_0}\,a_D &<& \sum_{i=1}^{l-1}{}'\left(\!\frac{\tilde\rho_{i+1}-
\tilde\rho_i}{k_i-k_{i+1}}-\frac{\tilde\rho_i-\tilde\rho_{i-1}}
{k_{i-1}-k_i}\right)\sum_{j=1}^{k_i}\big(N+1-h^0(\I^j_D(1))\big) \\
&=& \label{diverge} \sum_{i=1}^{l-1}{}'\ \frac{\tilde\rho_{i+1}-\tilde\rho_i}
{k_i-k_{i+1}}\sum_{j=k\_{i+1}+1}^{k_i}\big(N+1-h^0(\I^j_D(1))\big) \\ \nonumber
&\le& \sum_{i=1}^N(\tilde\rho_{i+1}-\tilde\rho_i)\big(N+1-h^0(\I_D^{k_i}(1))\big),
\end{eqnarray}
where in the last line we have added back in the $k_i=k_{i+1}$ terms since
they are positive.
Since the last sum is also $\sum_{i=2}^{N+1}\tilde\rho_i(h^0(\I_D^{k_i}(1))-
h^0(\I_D^{k_{i+1}}(1)))$, with a positive coefficient for each $\tilde\rho_i$,
we can replace $\tilde\rho_i$ by $\rho_i\ge\tilde\rho_i$ to give
$$
-\frac{N+1}{a_0}\,a_D<\sum_{i=1}^N(\rho_{i+1}-\rho_i)h^0(\O_{k_iD}(1)).
$$
Sum over the connected components $D$ and use $h^0(\I_i(1))\ge
i$ to give
\begin{eqnarray*}
-\frac{N+1}{a_0}\,a &<& \sum_{i=1}^N(\rho_{i+1}-\rho_i)h^0(\O_{Z_i}(1)) \\
&\le& \sum_{i=1}^N(\rho_{i+1}-\rho_i)\big(N+1-h^0(\I_i(1))+h^1(\I_i(1))\big) \\
&\le& \sum_{i=1}^N(\rho_{i+1}-\rho_i)\big(N+1-i\big)
+\sum_{i=1}^N(\rho_{i+1}-\rho_i)h^1(\I_i(1)) \\
&\le& \sum_{i=1}^{N+1}\rho_i+\delta.
\end{eqnarray*}
(In the above, any sum from $i$ to $j$ with $j<i$ is to be interpreted as
zero, and in passing from the first line to the second we have used the fact
that $\sum_{j=1}^{k_l}=0$ since $k_l=0$.)
\end{proof}

The Chow slope inequality (\ref{sss}) can be rewritten
$$
\epsilon+\frac{N+1}{a_0}=\epsilon+Ch(X)<Ch_c(\O_Z)=
\frac{\sum_{i=1}^c\big(N+1-h^0(\I_{\!Z}^i(1))\big)}{\int_0^c\!\tilde a_0(x)dx}
$$
for some small $\epsilon>0$. This implies
the weaker inequality (which is all that we shall require for curves)
\begin{equation} \label{css}
\left(\frac{N+1}{a_0}+\epsilon\right)\int_0^c\!\tilde a_0(x)dx<
\sum_{i=1}^ch^0(\O_{iZ}(1)).
\end{equation}
If $\epsilon$ can be chosen uniformly in (\ref{css}) for all $Z\subset X$
we call $X$ \emph{uniformly Chow slope stable with constant $\epsilon$}.
This will help us to deal with the correction $\delta$.

For $X$ a curve we can improve the estimates of Theorem \ref{uniform} slightly
and use Clifford's theorem to bound the $h^1$ terms:

\begin{thm} \label{uniformcurve}
If a smooth curve $(X,L)$ of genus $g$ and $d=\deg L>2g-2$ is
uniformly Chow slope stable (\ref{css}) with constant $\epsilon\ge
\big(1+\frac1{d-g}\big)\big(g-\frac12\big).\frac1d$ then it is Chow stable.
\end{thm}

\begin{proof}
We follow the same proof as above with $\frac{N+1}{a_0}$ replaced by
$\frac{N+1}{a_0}+\epsilon$ throughout, up to (\ref{diverge}), at which
point we use estimates specific to curves to better bound the $h^1$ terms.
That is, $D$ is locally a smooth point $\{p\}$ in $X$, and for those $i$
with $k_i\ne k_{i+1}$ (i.e.\ those involved in the sum $\sum{}'$),
$$
\frac1{k_i-k_{i+1}}\sum_{j=k\_{i+1}+1}^{k_i}h^0(\O_{j\{p\}}(1))=
\frac12\big(h^0(\O_{k_i\{p\}}(1))+h^0(\O_{k_{i+1}\{p\}}(1))+1\big).
$$
Therefore (\ref{diverge}) becomes
$$
-\left(\frac{N+1}{a_0}+\epsilon\right)a_p<\sum_{i=1}^{l-1}{}'
(\tilde\rho_{i+1}-\tilde\rho_i)\frac12\big(h^0(\O_{k_i\{p\}}(1))+
h^0(\O_{k_{i+1}\{p\}}(1))+1\big).
$$
We can now add back in those $i$ with $k_i=k_{i+1}$, as all terms are positive.
Each $\tilde\rho_i$ appears with positive coefficient in the result, since
it can be rearranged as $\sum_{i=2}^N\tilde\rho_i\big(
h^0(\O_{k_{i-1}\{p\}}(1))-h^0(\O_{k_{i+1}\{p\}}(1))\big)+\rho\_{\!N+1}\big(
h^0(\O_{k_N\{p\}}(1))+1\big)$.
So replacing $\tilde\rho_i$ by $\rho_i\ge\tilde\rho_i$, summing over
$p$ in the support
of $Z_1$, and using $h^0(\I_i(1))\ge i$, gives
\enlargethispage*{1ex}
\begin{eqnarray*}
-\left(\frac{N+1}{a_0}+\epsilon\right)a \hspace{-2cm} & \\ \qquad\quad &<&
\sum_{i=1}^N(\rho_{i+1}-\rho_i)\frac12
\big(h^0(\O_{Z_i}(1))+h^0(\O_{Z_{i+1}}(1))+1\big) \\
&\le& \sum_{i=1}^N(\rho_{i+1}-\rho_i)\frac12\Big[N+1-h^0(\I_i(1))
+h^1(\I_i(1))+ \\
&& \hspace{45mm} N+1-h^0(\I_{i+1}(1))+h^1(\I_{i+1}(1))+1\Big] \\
&\le& \sum_{i=1}^N(\rho_{i+1}-\rho_i)\big[N+1-i\big]
+\frac12\sum_{i=1}^N(\rho_{i+1}-\rho_i)\big(h^1(\I_i(1))+h^1(\I_{i+1}(1))\big)
\\ &=& \sum_{i=1}^{N+1}\rho_i+\delta, \hspace{20mm}
\delta:=\frac12\sum_{i=1}^N(\rho_{i+1}-\rho_i)
\big(h^1(\I_i(1))+h^1(\I_{i+1}(1))\big)
\end{eqnarray*}
i.e.
\begin{equation}
-\left(1+\epsilon\frac{a_0}{N+1}\right)\frac{N+1}{a_0}\,a
<\left(1+\frac{\delta}{\sum_{i=1}^{N+1}\rho_i}\right)
\sum_{i=1}^{N+1}\rho_i\,. \label{d}
\end{equation}
This $\delta$ is a tiny improvement over the one in Theorem \ref{uniform},
but can be bounded using Clifford's theorem, which for our purposes says
that $h^1(L)\le\max\{1+g-h^0(L),0\}$. Thus $h^1(\I_i(1))\le1+g-i$ for $i\le
g$ and vanishes for $i>g$. This yields
\begin{multline*}
\delta\le\sum_{i=1}^g(\rho_{i+1}-\rho_i)\frac12\big(1+g-i+1+g-(i+1)\big)=
\\ \sum_{i=1}^g(\rho_{i+1}-\rho_i)\Big(g-i+\frac12\Big)=\sum_{i=2}^g\rho_i+
\frac12\rho_{g+1},
\end{multline*}
since $\rho_1=0$. Since the $\rho_i$ are monotonic in $i$,
$$
\sum_{i=2}^g\rho_i+\frac12\rho_{g+1}\le\Big(g-\frac12\Big)\rho_{g+1}
\le\Big(g-\frac12\Big)\frac{\frac12\rho_{g+1}+\sum_{i=g+2}^{N+1}\rho_i}
{N-g+\frac12}\,.
$$
Adding $\frac{g-\frac12}{N-g+\frac12}\big(\sum_{i=2}^g\rho_i+
\frac12\rho_{g+1}\big)$ to both sides gives
$$
\frac N{N-g+\frac12}\left(\sum_{i=2}^g\rho_i+\frac12\rho_{g+1}\right)
\le\frac{g-\frac12}{N-g+\frac12}\sum_{i=2}^{N+1}\rho_i,
$$
and so
$$
\delta\le\sum_{i=2}^g\rho_i+\frac12\rho_{g+1}\le
\frac{g-\frac12}N\sum_{i=2}^{N+1}\rho_i.
$$
Combined with (\ref{d}) we find that uniform slope stability implies that
$$
-\left(1+\epsilon\frac{a_0}{N+1}\right)\frac{N+1}{a_0}\,a<
\left(1+\frac{g-\frac12}N\right)\sum_{i=1}^{N+1}\rho_i,
$$
which implies the inequality $-\frac{N+1}{a_0}\,a<\sum_{i=1}^{N+1}\!\rho_i\,$
required (\ref{mummy}) if $\epsilon\ge\frac{N+1}{a_0}\frac{g-\frac12}N$.
The condition $d>2g-2$ implies that $h^1(\O(1))=0$ so $N+1=d+1-g$ and $a_0=d$.
Therefore the inequality is $\epsilon\ge\big(1+\frac1{d-g}\big)
\big(g-\frac12\big)\frac1d$.
\end{proof}

\begin{thm}
Smooth curves $(X,\O(1))$ of genus $g\ge1$ are uniformly slope stable (\ref{css})
with any constant $\epsilon<\frac gd$, and so are asymptotically Chow stable.
\end{thm}

\begin{proof}
We need only demonstrate the inequality
(\ref{css}) in the case of $Z=\{p\}$ a single reduced point, as the inequality
for a multiple of $\{p\}$ is weaker and the inequality for arbitrary $Z$
follows from adding the inequalities for its connected components.
By Riemann-Roch, $\tilde a_0(x)=x$, so that $\int_0^c\tilde a_0(x)=c^2/2$
while $\sum_{i=1}^ch^0(\O_{iZ}(1))=\sum_{i=1}^ci=c(c+1)/2$. Therefore
$$
\left(\frac{d+1-g}d+\epsilon\right)\frac{c^2}2=
\left(\frac{N+1}{a_0}+\epsilon\right)\int_0^c\!\tilde a_0(x)dx\ <\ 
\sum_{i=1}^ch^0(\O_{iZ}(1))=\frac{c^2+c}2\,,
$$
so long as $\epsilon<\frac1c+\frac{g-1}d$. Of course $c\le d$, so $(X,L)$
is uniformly Chow slope stable with any constant $\epsilon<\frac gd$.

Theorem \ref{uniformcurve} then gives Chow stability, so long as
$\frac gd>\big(1+\frac1{d-g}\big)\big(g-\frac12\big).\frac1d$ which is true
for $g\ge1$ and sufficiently large $d$.
\end{proof}

\begin{rmk}
  Mumford \cite{Mu} proves the sharper result that $(X,\O(1))$ is Chow stable
  for $\deg\O(1)>2g$ using a combinatorial argument.
\end{rmk}

\section{Examples} \label{egs}

Our remaining examples all deal with K-slope stability. Many more examples,
calculations and applications are given in \cite{RT}.

\subsection{Varieties with nonnegative canonical bundle}\ \vskip 5pt
\noindent
Suppose that $X$ has at worst canonical singularities.
That is $X$ is normal, there is an integer $m$ such that $mK_X$ is
Cartier, and given any resolution of singularities
$\pi_1\colon\overline{\!X}\,\rightarrow X$ we have
\begin{equation}\label{canonical}
mK_{\overline{\!X}}=\pi_1^*(mK_X)+\sum_i\alpha_iF_i \quad\text{with}\ \ \alpha_i\ge 0, 
\end{equation}
where the $F_i$ are the irreducible components of the
exceptional set of $\pi_1$. We can define intersection with the canonical
class of $X$ on $\overline{\!X}\,$ by $K_X.(\,\cdot\,):=\frac1m(mK_X).(\,\cdot\,)$.

For any subscheme $Z\subset X$ let
$$
\overline{\!X}\Rt{\pi_2\ }\widehat X\rt{\pi}X
$$
be a resolution of singularities of $\widehat X$, the blow up of $X$ along
$Z$. $\pi_1=\pi\comp\pi_2\colon\overline{\!X}\rightarrow
X$ is a resolution of singularities of $X$, so (\ref{canonical}) holds.
Letting $F=\pi_2^*E$, $L-xF=\pi_2^*(L-xE)$ is nef on $\overline{\!X}$
for $0\le x\le\epsilon(Z)$.

We wish to compute $a_i(x)$ on $\overline{\!X}$ instead of $\widehat
X$. Since $\O_{\widehat X}\subseteq\pi_{2*}\O_{\,\overline{\!X}}$ with quotient
supported in codimension one, we have an inclusion
$$
H^0_{\widehat X}((L-xE)^k)\subseteq H^0_{\overline{\!X}}((L-xF)^k),
$$
with cokernel of dimension $\le O(k^{n-1})$ by Fujita vanishing, for
$0\le x<\epsilon(Z)$. As in (\ref{normal}), for $k\gg0$, $h_{\,\overline{\!X}}^i
((L-xF)^k)=h^0_{\widehat X}((L-xE)^k\otimes R^i\pi_{2*}\O)=O(k^{n-1-i})$
since the support of $R^i\pi_{2*}\O$ has codimension at least $i+1$.
Therefore $\chi\_{\widehat X}((L-xE)^k)=\chi\_{\overline{\!X}}((L-xF)^k)
-ak^{n-1}+O(k^{n-2})$ for some $a\ge0$. That is, for $0\le x<\epsilon(Z)$,
we have
$$
a_0(x) = \frac{1}{n!}(L-xF)^n, \vspace{-6pt}
$$
and \vspace{-6pt}
\begin{equation}\label{a_idefinedbyaresolution}
\qquad\qquad a_1(x)\le-\frac{1}{2(n-1)!}K_{\overline{\!X}}.(L-xF)^{n-1}
\le-\frac{1}{2(n-1)!}K_{X}.(L-xF)^{n-1},
\end{equation}
where the second inequality follows from (\ref{canonical}) and the
fact that $L-xF$ is nef. Again as in (\ref{normal}), the first two
terms of the Euler characteristic of $L$ are the same as those of
$\pi_1^* L$, so equality holds in (\ref{a_idefinedbyaresolution}) when
$x=0$.  That is $a_0=\frac1{n!}L^n$ and
\begin{equation} \label{a_1definedbyaresolution}
a_1 = -\frac{1}{2(n-1)!}  K_{\overline{\!X}}.L^{n-1}=- \frac{1}{2(n-1)!} K_X.L^{n-1},
\end{equation}
since $L$ is trivial along the $F_i$. With these preliminaries, it becomes
easy to prove the following.

\begin{thm} \label{CY} \textbf{Calabi-Yaus and canonical models}. \\
Let $X$ be an irreducible variety with at worst canonical singularities.
\newline
${}\quad\bullet$ If $K_X$ is numerically trivial then $(X,L)$ is slope stable for all polarisations $L$. \newline
${}\quad\bullet$ If $K_X$ is ample then the canonical polarisation $(X,K_X)$
is slope stable.
\end{thm}

\begin{proof}
In both cases, $K_X\sim\alpha L$ is numerically equivalent to a nonnegative
multiple $\alpha\ge0$ of the polarisation. So by (\ref{a_1definedbyaresolution}),
$\mu(X)=a_1/a_0=-n\alpha/2$. By (\ref{a_idefinedbyaresolution}), 
$$
-\mu(X)a_0(x)+a_1(x)\le\frac{\alpha}{2(n-1)!}(L-xF)^n-\frac1{2(n-1)!}
(\alpha L).(L-xF)^{n-1},
$$
which equals $-\frac{\alpha x}{2(n-1)!}F.(L-xF)^{n-1}\le0$ since $L-xF$ is
nef. Since $a_0'(x)<0$ for $x\in(0,\epsilon(Z))$ (\ref{a0'}), integration
gives
$$
-\mu(X)\int_0^ca_0(x)dx+\int_0^ca_1(x)+\frac{a_0'(x)}2dx<0 \quad\text{for}\
c\in(0,\epsilon(Z)],
$$
which rearranges to give slope stability: $\mu_c(\I_Z)<\mu(X)$.
\end{proof}

With more work, this can be generalised as follows.

\begin{thm}
Suppose that $(X,L)$ has at worst canonical singularities, and
$K_X$ is nef and big. Then $(X,L)$ is slope stable
for $L$ ample and sufficiently close to $K$.  More precisely,
  \begin{itemize}
  \item For any divisor $G$ there is a $\delta_0>0$ such
    that if $0\le \delta<\delta_0$ and $L=K_X+\delta G$ is ample then
    $(X,L)$ is slope stable.
  \item If $2\mu(X,L)L + nK_X$ is nef then $(X,L)$ is slope stable.
  \item If $-2\mu(X,L)L - nK$ is nef then $(X,L)$ is slope stable.
\end{itemize}
\end{thm}

In \cite{RT} we prove that no smooth $Z$ can slope destabilise a smooth $(X,L)$
with these properties, and the proof extends to general $Z$ and $X$ with
canonical singularities using the preliminaries
(\ref{a_idefinedbyaresolution}) and (\ref{a_1definedbyaresolution}).

These results are to be expected due to the deep and difficult related
results of Viehweg \cite{V}, and the expectation that the minimal
model programme can be carried out in all dimensions (and so can be
done in families) \cite{Ka}.  In fact Theorem \ref{CY} can be proved in a
round about way for smooth varieties with no holomorphic vector fields by
the stability results of \cite{Do1, Zh} applied to the K\"ahler-Einstein
metrics of \cite{Au, Y} on such
varieties. Similarly, if $X$ has no holomorphic vector fields and $L$ is
ample and sufficiently close to $K_X$ then an implicit function theorem
argument applied to the K\"ahler-Einstein metric provides a constant scalar
curvature K\"ahler metric in $[c_1(L)]$.  But the quick proofs
above demonstrate that using slope stability it could become much easier
to produce and compactify moduli of varieties with semi-ample canonical
bundle \cite{V}.

\subsection{Irreducible Curves} \label{curves}\ \vskip 5pt

\begin{prop}\label{prop:stabilityofpolarisedcurves}
  Let $(\Sigma,L)$ be an irreducible polarised curve with arithmetic genus
  $g\ge 2$. The Hilbert-Samuel polynomial of a subscheme $Z\subset\Sigma$
  can be written  
\begin{equation}\label{hilbertsamuelofsubschemesofcurves} 
h^0(L^k\big/(L^k\otimes\I_{\!Z}^{xk})) = e(Z)xk-\rho(Z)\quad \text{for}\ k\gg0,
\end{equation} 
  and $Z$ destabilises $\Sigma$ if and only if it strictly
  destabilises, if and only if $2\rho(Z)> e(Z)$.
\end{prop}

\begin{proof}
  The assumption on the genus implies that $\mu(\Sigma)<0$. As $\dim 
\Sigma=1$, 
$\tilde{a}_0(x)$
  is a degree 1 polynomial vanishing at the origin, while $\tilde a_1(x)$
  has degree 0. So writing $\tilde a_0(x)=e(Z)x$ and $\tilde{a}_1(x)=-\rho(Z)$,
  $$
  \mu_c(\O_Z) = \frac{\int_0^c \tilde{a}_1(x) +
    \frac{\tilde{a}_0'(x)}{2} dx}{\int_0^c \tilde{a}_0(x) dx} =
  \frac{e(Z)-2\rho(Z)}{e(Z)c}\,.
  $$
  This has the same sign as $e(Z)-2\rho(Z)$ for all $c>0$, and tends to $\pm\infty$
  as $c\to0$. Thus $p$ destabilises if and only if it strictly destabilises
  if and only if this sign is negative.
\end{proof}

\begin{thm} \label{multiplicityofdestabilisingsingularpoints}
  Let $(\Sigma,L)$ be a polarised irreducible curve of arithmetic
  genus $g\ge 2$.
  \begin{itemize}
  \item If $\Sigma$ has a point of multiplicity $e\ge 3$ then
    $(\Sigma,L)$ is not slope stable.
  \item If $\Sigma$ has at worst ordinary double points then
    $(\Sigma,L)$ is slope stable.
\end{itemize}
\end{thm}

\begin{proof}
  For the first statement suppose $Z=\{p\}$ is a point of $\Sigma$
  with multiplicity $e\ge 3$.  Northcott has shown (\cite{No}
  Lemma 1) that $\rho\ge e-1$.  Hence $2\rho\ge 2e-2>e$ and
  $Z$ strictly destabilises.
  
  For the second statement suppose for a contradiction that $Z$ is a
  destabilising subscheme of $\Sigma$.  From Proposition
  \ref{simplify} (especially (\ref{simplify2})) and Lemma
  \ref{seshadridisconnected} we see that
  $\mu_c(\I_{Z_0})\ge\mu(\Sigma)$ for some connected component $Z_0$ of $Z$
  and $c\le\epsilon(Z_0)$ (we are not saying that $Z_0$ destabilises since
  we may not be allowed to take $c=\epsilon(Z_0)$). If the support $p\in\Sigma$
  of $Z_0$ is a smooth point then $Z_0=m\{p\}$ must by a thickened point,
  which
  by Proposition \ref{simplify} shows that $\mu_{mc}(\I_{\{p\}})\ge\mu(X)$.
  But in the notation of (\ref{hilbertsamuelofsubschemesofcurves}),
  $e(\{p\})=1$, which by Theorem 3.2 of \cite{KM} implies that
  $\rho(\{p\})=0$, so by Proposition \ref{prop:stabilityofpolarisedcurves}
  $\mu_{mc}(\I_{\{p\}})$ is in fact $<\mu(X)$.
  
  Hence $Z$ must be supported at one of the singular points, and we
  can reduce to looking at the local analytic model $\Spec R,\ 
  R=\C[X,Y]/(XY)$.  Let $I$ be an ideal of $R$
  which is supported at $(X,Y)$.  (By abuse of notation we shall not
  distinguish between a polynomial in two variables and its class in
  $R$.)  We can pick a finite number of generators of $I$ of the form
  $f_i=a_iX^{p_i}+b_iY^{q_i}$, with $a_i,b_i\in\mathbb C$.
   
  Let $p=\min\{p_j\colon a_j\neq 0\},\ q=\min\{q_j\colon a_j\neq 0\}$;$p,q\ge1$ since $I$ is supported at the origin.
  Pick an $i$ such that $p_i=p$ and $a_i\ne0$;
  then $X^{pk+1}=\frac{X}{a_i^k}(a_iX^{p_i} + b_i Y^{q_i})^k\in I^k$; similarly
  $Y^{qk+1}\in I^k$. Therefore
  $R/I^k$ is spanned by $\{1,X,\ldots,X^{pk},Y,\ldots
  Y^{qk}\}$.  By the definition of $p$ and $q$, $I^k$ is spanned by
  $\{X^i, Y^j \colon i\ge pk, j\ge qk\}$, so the vectors $\{1,X,\ldots
  X^{pk-1}, Y,\ldots, Y^{qk-1}\}$ in $R/I^k$ are linearly independent.
  Thus
  $$
  (p+q)k-1\le \dim R/I^k \le (p+q)k+1.
  $$
    Writing $\dim R/I^k=ek-\rho$ we have $-1\le\rho\le 1$.
  Hence $2\rho \le 2\le p+q=e$, so by Proposition
  \ref{prop:stabilityofpolarisedcurves} $I$ does not destabilise.
\end{proof}

\begin{rmk}
  Eisenbud and Mumford \cite{Mu} analyse the effect of singular points
  on Chow stability of higher dimensional varieties. It would be interesting
  to know if their results can be seen using slope stability.
\end{rmk}

These results combined with Theorem \ref{thm:kstableslopestable} imply
that curves with singularities of multiplicity greater than two are strictly
K-unstable. We cannot deduce positive results about K-stability from the
results of Section \ref{converse}, however, since $\Sigma$ need not
be normal. Unless, that is, $\Sigma$ is smooth:

\begin{thm}\label{smoothcurvesareKstable} Any smooth polarised curve $(\Sigma,L)$
of genus $g$ is K-stable if $g\ge1$ and strictly K-polystable if $g=0$.
\end{thm}

\begin{proof}
By Corollary \ref{Kslopecurve} it is equivalent to prove the results for
slope (poly)stability. Instead of using previous results it is now easier
to proceed directly. Any nonempty subscheme $Z$ is a divisor
of degree $d>0$, so
$$
\chi(L^k\otimes\I_{\!Z}^{xk})=k\deg L-xdk+1-g
$$
shows that $\tilde a_0(x)=xd$ and $\tilde a_1(x)=0$. Thus $\mu_c(\O_Z)=
\frac{cd}{c^2d}=\frac1{c}>0\ge\frac{1-g}{\deg L}=\mu(X)$ for $g\ge1$, proving
slope stability.

For $g=0$, $c$ may take values up to and \emph{including} $\epsilon(Z)=\deg
L/d$, since $L^d\otimes\I_Z^{\deg L}=\O_{\PP^1}(d\deg L-d\deg L)=\O_{\PP^1}$
is globally generated.
Thus $\mu_c(\O_Z)\ge\frac{d}{\deg L}\ge\frac1{\deg L}=\mu(X)$ with equality
(strict semistability) only for $d=1$, \emph{i.e.}\ $Z$ a single point, and
$c=\epsilon(Z)$. Since the deformation to
  the normal cone of a single point on $\PP^1$ blows down to $\PP^1\times\C$
  (with a nontrivial $\C^\times$-action) from which the relevant line bundle
  $\L_c$ pulls back, we find $\PP^1$ is in fact slope polystable.
\end{proof} 

This can also be proved using the constant curvature
metric on $\Sigma$ and analysis of the Mabuchi functional, but this
seems to be the first direct algebraic proof.

\vskip 4mm

{\small \noindent {\tt jaross@math.columbia.edu} \\
\noindent \small{\tt richard.thomas@imperial.ac.uk}} \newline
\noindent Department of Mathematics, Columbia University, New York, NY 10027.
USA. \\
Department of Mathematics, Imperial College, London SW7 2AZ. UK.

\end{document}